\documentclass{Article2022}

\newcommand{\Operad}{\mathcal{O}}
\DeclareMathOperator{\Composition}{\circ}
\newcommand{\SetCliff}{\mathsf{Cl}}
\newcommand{\SetHill}{\mathsf{Hi}}
\newcommand{\SetPrime}{\mathcal{P}}
\newcommand{\Leq}{\preccurlyeq}
\newcommand{\JJoin}{\vee}
\newcommand{\Meet}{\wedge}
\newcommand{\GeneratingSet}{\mathfrak{G}}
\newcommand{\RelationSpace}{\mathcal{R}}
\newcommand{\FussCatalan}{\mathrm{cat}}
\newcommand{\Reduction}{\mathrm{r}}
\newcommand{\SpaceCliff}{\mathbf{Cl}}
\newcommand{\SpaceHill}{\mathbf{Hi}}
\newcommand{\BasisE}{\mathsf{E}}
\newcommand{\BasisF}{\mathsf{F}}
\newcommand{\BasisH}{\mathsf{H}}
\newcommand{\BasisK}{\mathsf{K}}
\newcommand{\SpaceV}{\mathcal{V}}
\newcommand{\SubFamilly}{\mathcal{S}}
\newcommand{\CharacteristicCliff}{\chi}
\newcommand{\IntersticeOperad}{\mathbf{I}}
\newcommand{\Length}{\ell}
\newcommand{\MobiusFunction}{\mu}
\newcommand{\As}{\mathbf{As}}
\newcommand{\ConstantLength}{\mathrm{c}}
\newcommand{\DominationLength}{\mathrm{d}}
\newcommand{\Dual}{\star}
\newcommand{\MapConstant}[1]{{\underline{#1}}}
\newcommand{\MapArithmetic}[1]{{\mathbf{#1}}}
\newcommand{\Modify}{\mathrm{m}}
\newcommand{\LLambda}{\bm{\lambda}}
\newcommand{\ToComposition}{\mathrm{comp}}
\newcommand{\ToPermutation}{\mathrm{perm}}
\newcommand{\ToTree}{\mathrm{tree}}
\newcommand{\ToPath}{\mathrm{path}}
\newcommand{\ToDyck}{\mathrm{dyck}}
\newcommand{\TreeT}{\mathfrak{t}}

\newcommand{\Span}{\mathrm{Span}}
\newcommand{\Unit}{\mathds{1}}
\newcommand{\GenA}{\mathtt{a}}
\newcommand{\GenB}{\mathtt{b}}
\newcommand{\GenC}{\mathtt{c}}

\DeclareMathOperator{\WhiteSquare}{
\begin{tikzpicture}[Centering,scale=.16]
    \node[rectangle,draw=ColA!55,fill=ColA!5,thick,inner sep=0,minimum size=1.75mm](0,0){};
\end{tikzpicture}}

\DeclareMathOperator{\BlackSquare}{
\begin{tikzpicture}[Centering,scale=.16]
    \node[rectangle,draw=ColB!100,fill=ColB!45,thick,inner sep=0,minimum size=1.75mm]
        (0,0){};
\end{tikzpicture}}

\title[Cliff operads]{Cliff operads: a hierarchy of operads on words \vspace{-5ex}}
\date{\today}
\author[C. Combe]{Camille Combe}
\address{
    \vspace{-0.4cm}
    \scriptsize
    LIGM, Université Gustave Eiffel, CNRS, ESIEE Paris, F-$77454$
    Marne-la-Vallée, France.
    {\tt \href{mailto:camille.combe@univ-eiffel.fr}{camille.combe@univ-eiffel.fr}}}
\author[S. Giraudo]{\vspace{-0.8cm} Samuele Giraudo}
\address{
    \vspace{-0.4cm}
    \scriptsize
    LIGM, Université Gustave Eiffel, CNRS, ESIEE Paris, F-$77454$
    Marne-la-Vallée, France.
    {\tt \href{mailto:samuele.giraudo@univ-eiffel.fr}{samuele.giraudo@univ-eiffel.fr}}}
\thanks{%
    This research has been partially supported by the projects CARPLO (ANR-20-CE40-0007)
    and ALCOHOL (ANR-19-CE40-0006) of the Agence nationale de la recherche.}

\begin{document}

\begin{abstract}
    A new hierarchy of operads over the linear spans of $\delta$-cliffs, which are some
    words of integers, is introduced. These operads are intended to be analogues of the
    operad of permutations, also known as the associative symmetric operad. We obtain
    operads whose partial compositions can be described in terms of intervals of the lattice
    of $\delta$-cliffs. These operads are very peculiar in the world of the combinatorial
    operads since, despite the relative simplicity for their construction, they are
    infinitely generated and they have nonquadratic and nonhomogeneous nontrivial relations.
    We provide a general construction for some of their quotients. We use it to endow the
    spaces of permutations, $m$-increasing trees, $c$-rectangular paths, and $m$-Dyck paths
    with operad structures. The operads on $c$-rectangular paths admit, as Koszul duals,
    operads generalizing the duplicial and triplicial operads.
    \vspace{-4ex}
\end{abstract}

\MakeFirstPage

\section*{Introduction}
Endowing sets of combinatorial objects with operations has a long-term tradition. The
literature abounds of examples of monoids, lattices, pre-Lie algebras, associative algebras,
Hopf bialgebras, and operads defined on the linear span of combinatorial sets (see among
others the recent works~\cite{MR95,LR02,BF03,AL07,Gir12,Gir15,LPRM20}). Adding such an
algebraic dimension offers a new point of view of the objects: we can see these as
assemblies, through the offered operations, of generators of the considered structure. This
gives among others tools to enumerate combinatorial sets or to establish transformations
(and in particular bijections) between two sets of combinatorial objects. All this also
maintains connections with partial order theory because the operations in most of these
structures can be described as intervals of some partial orders on the underlying objects
(as few examples, see~\cite{LPRM20} for some operads on $m$-Dyck paths, and~\cite{LR02} for
Hopf bialgebras on permutations and binary trees). In this vein, there is a very rich operad
structure on the linear span $\As$ of all permutations for which the operadic partial
composition can be described as intervals of the right weak order~\cite{AL07}.

The primary impetus for this work was the will to introduce a similar operad structure on
the linear span of the set of all $\MapArithmetic{1}$-cliffs, a set of words of integers
which is in one-to-one correspondence with the set of permutations. These objects admit some
generalizations named $\delta$-cliffs, depending on a parameter $\delta$ which is a map from
$\N \setminus \{0\}$ to $\N$ (or equivalently, an infinite word of integers). In our
context, we search for an analogue of $\As$ involving $\delta$-cliffs instead, and for which
the operadic partial composition can be described as intervals of a lattice on
$\delta$-cliffs introduced in~\cite{CG20,CG22}. As shown in this previous work of the two
present authors, the linear span $\SpaceCliff_\delta$ of all $\delta$-cliffs when $\delta$
is unimodal admits the structure of an associative algebra. In the present work, we show
that this space admits (up to a shift in its graduation) also the structure of an operad.
Surprisingly, the construction works again only when $\delta$ is unimodal. In a remarkable
way, in the same way as for $\As$, the partial composition of $\SpaceCliff_\delta$ can be
described in terms of intervals of the poset of $\delta$-cliffs. More precisely, there is a
basis of $\SpaceCliff_\delta$ for which the partial composition of two basis elements admits
as support the empty set or an interval of the poset of $\delta$-cliffs.

The main results presented here include the construction of three different bases of
$\SpaceCliff_\delta$, a necessary and sufficient condition on $\delta$ for the fact that
$\SpaceCliff_\delta$ is finitely generated, and a sufficient condition on $\delta$ for the
fact that the space of nontrivial relations of $\SpaceCliff_\delta$ is not finitely
generated. We also explore a way to construct quotient operads $\SpaceCliff_\SubFamilly$ of
$\SpaceCliff_\delta$ whose bases are indexed by particular subsets $\SubFamilly$ of
$\delta$-cliffs. We show here that when $\SubFamilly$ is a sublattice of the lattice of
$\delta$-cliffs, the partial composition of $\SpaceCliff_\SubFamilly$ can be described in
terms of intervals of $\SubFamilly$. We finally explore some concrete examples of operads
$\SpaceCliff_\delta$. These operads appear, unexpectedly, to have a very rich and complex
structure. For instance, the space of nontrivial relations of
$\SpaceCliff_\MapArithmetic{1}$ is infinitely generated and has elements which are
nonquadratic and nonhomogeneous in terms of degrees. This is quite rare in the panorama of
combinatorial operads and seems to brought to the light a very singular new object. The
bases of the operads $\SpaceCliff_\MapArithmetic{m}$, where $\MapArithmetic{m}$ is the
arithmetic sequence with a common difference of $m$, are index by some labeled planar rooted
trees and have also an intricate structure. We also explore the case of the quotients
$\SpaceCliff_\SubFamilly$ where $\SubFamilly$ is the set of weakly increasing
$\delta$-cliffs, called $\delta$-hills in~\cite{CG20,CG22}. We obtain here operads whose
dimensions are provided by shifted Fuss-Catalan numbers (which are hence different from the
operads of~\cite{LPRM20} whose dimensions are given by Fuss-Catalan numbers) and where the
Stanley lattice~\cite{Sta75} is the underlying partial order for the description of the
partial composition. We also construct a last family of operads whose dimensions are
provided by some binomial coefficients and whose bases are indexed by some paths formed by
east and north steps in rectangles. They have the interesting property to have, as Koszul
duals, operads having dimensions enumerated by (not shifted) Fuss-Catalan numbers which can
be thought as generalizations of the duplicial operad~\cite{BF03} and the triplicial
operad~\cite{Ler11}.

This paper is organized as follows. Section~\ref{sec:preliminaries} exposes some definitions
about $\delta$-cliffs and some background and notations about operads.  In
Section~\ref{sec:cliff_operads}, we provide the construction of the operads
$\SpaceCliff_\delta$ and their first properties. Section~\ref{sec:quotient_operads} is
devoted to the study of the quotient operads $\SpaceCliff_\SubFamilly$ of
$\SpaceCliff_\delta$. Section~\ref{sec:particular_cases} presents some particular examples
arising from our constructions.

\subsubsection*{General notations and conventions}
For any integers $i$ and $j$, $[i, j]$ denotes the set $\{i, i + 1, \dots, j\}$. For any
integer $i$, $[i]$ denotes the set $[1, i]$ and $\HanL{i}$ denotes the set $[0, i]$. For any
set $A$, $A^*$ is the set of all words on $A$. If $a$ is a letter and $n$ is a nonnegative
integer, $a^n$ is the word consisting in $n$ occurrences of $a$. In particular, $a^0$ is the
empty word~$\epsilon$. For any $w \in A^*$, $\Length(w)$ is the length of $w$, and for any
$i \in [\Length(w)]$, $w(i)$ is the $i$-th letter of $w$. For any $i \leq j \in
[\Length(w)]$, $w(i, j)$ is the word $w(i) w(i + 1) \dots w(j)$. All algebraic structures
considered in this work have a field $\K$ of characteristic zero as ground field.

\section{Preliminaries} \label{sec:preliminaries}
This first part contains elementary definitions about cliffs and related combinatorial
objects. We also provide brief recalls about nonsymmetric operads and finish by describing
interstice operads. These operads contain, as suboperads or quotients, the forthcoming
operads on cliffs.

\subsection{Cliffs and related objects} \label{subsec:cliffs_and_objects}
Cliffs are essentially words of nonnegative integers satisfying some
properties~\cite{CG20,CG22}. We give here basic definitions about these objects and explain
how particular families of cliffs can encode other known combinatorial families (like
integer compositions, permutations, $m$-increasing trees, $c$-rectangular paths, and
$m$-Dyck paths).

\subsubsection{Graded sets}
A \Def{graded set} is a set expressed as a disjoint union
\begin{equation}
    S := \bigsqcup_{n \in \N} S(n)
\end{equation}
such that all $S(n)$, $n \in \N$, are sets. The \Def{size} $|x|$ of an $x \in S$ is the
unique integer $n$ such that $x \in S(n)$. If $S$ and $S'$ are two graded sets, a map
$\theta : S \to S'$ is a \Def{graded set morphism} if for any $x \in S$, $|\theta(x)| =
|x|$. Besides $S'$ is a \Def{graded subset} of $S$ if for any $n \in \N$, $S'(n) \subseteq
S(n)$.

\subsubsection{Cliffs}
A \Def{range map} is a map $\delta : \N \setminus \{0\} \to \N$.  For any $m \geq 0$, we
denote by $\MapArithmetic{m}$ the range map satisfying $\MapArithmetic{m}(i) = (i - 1)m$ for
any $i \geq 1$.  Moreover, for any $c \geq 0$, we denote by $\MapConstant{c}$ the range map
satisfying $\MapConstant{c}(i) = c$ for any $i \geq 1$. We shall specify range maps as
infinite words $\delta = \delta(1) \delta(2) \dots$. For this purpose, for any $a \in \N$,
we shall denote by $a^\omega$ the infinite word having all its letters equal to $a$. For
instance, the notation $\delta := 2113^\omega$ stands the range map $\delta$ satisfying
$\delta(1) = 2$, $\delta(2) = \delta(3) = 1$, and $\delta(i) = 3$ for all $i \geq 4$. A
range map $\delta$ is \Def{unimodal} if for any $1 \leq i_1 < i_2 < i_3$, $\delta\Par{i_1} >
\delta\Par{i_2} < \delta\Par{i_3}$ never occurs. Besides, $\delta$ is \Def{$1$-dominated} if
there is a $k \geq 1$ such that for all $k' \geq k$, $\delta(1) \geq \delta\Par{k'}$.

A word $w$ on $\N$ is a \Def{$\delta$-cliff} if for any $i \in [\Length(w)]$, $w(i) \in
\HanL{\delta(i)}$. The \Def{size} $|w|$ of a $\delta$-cliff $w$ is $\Length(w) + 1$. The
graded set of all $\delta$-cliffs is denoted by $\SetCliff_\delta$. For instance
$\SetCliff_\MapArithmetic{1}(4) = \{000, 001, 002, 010, 011, 012\}$. Let also
$\SetHill_\delta$ be the graded subset of $\SetCliff_\delta$ containing all weakly
increasing $\delta$-cliffs. For instance $\SetHill_\MapArithmetic{1}(4) = \{000, 001, 002,
011, 012\}$. Any element of $\SetHill_\delta$ is a \Def{$\delta$-hill}. The
\Def{$\delta$-reduction} of a word $w$ on $\N$ is the $\delta$-cliff $\Reduction_\delta(w)$
satisfying
\begin{equation}
    \Par{\Reduction_\delta(w)}(i) = \min \Bra{w(i), \delta(i)}
\end{equation}
for any $i \in [\Length(w)]$. For instance, $\Reduction_\MapArithmetic{1}({\bf 2} 120 {\bf
6} {\bf 6}) = 012045$ and $\Reduction_\MapArithmetic{2}({\bf 2} 12066) = 012066$.

Let $\Leq$ be the partial order relation on $\SetCliff_\delta$ satisfying $u \Leq v$ for any
$u, v \in \SetCliff_\delta$ such that $|u| = |v|$ and for all $i \in [\Length(u)]$, $u_i
\leq v_i$. For any $n \geq 1$, the poset $\Par{\SetCliff_\delta(n), \Leq}$ is the
\Def{$\delta$-cliff poset} of order~$n$. This poset is a Cartesian product of total orders.
For this reason, it is a distributive lattice. For any $u, v \in \SetCliff(n)$ such that $u
\Leq v$, we denote by $[u, v]_\Leq$ the interval between $u$ and $v$ in this poset.

We now present graded sets which are in one-to-one correspondence with some particular sets
of cliffs.

\subsubsection{Integer compositions} \label{subsubsec:integer_compositions}
For any $n \geq 1$, $\SetCliff_\MapConstant{1}(n)$ is the set of all binary words of length
$n - 1$. For this reason, this set is in one-to-one correspondence with the set of all
\Def{integer compositions} of $n$, which are sequences $\Par{\LLambda_1, \dots, \LLambda_k}$
of positive integers such that $\LLambda_1 + \dots + \LLambda_k = n$. A possible bijection
sends $w \in \SetCliff_\MapConstant{1}(n)$ to the integer composition $\ToComposition(w) :=
\Par{\LLambda_1, \dots, \LLambda_k}$ where $w$ has $k - 1$ occurrences of $1$ and provided
that $w = 0^{\LLambda_1 - 1} 1 0^{\LLambda_2 - 1} 1 \dots 1 0^{\LLambda_k - 1}$. For
instance, in $\SetCliff_\MapConstant{1}(8)$,
\begin{equation}
    \ToComposition(1100010) = (1, 1, 4, 2).
\end{equation}
Observe that the poset $\Par{\SetCliff_\MapConstant{1}(n), \Leq}$ is the Boolean lattice of
order~$n - 1$.

\subsubsection{Permutations} \label{subsubsec:permutations}
For any $n \geq 1$, $\SetCliff_\MapArithmetic{1}(n)$ is the set of all words $w$ of length
$n - 1$ such that $w(i) \in \HanL{i - 1}$ for all $i \in [n - 1]$. For this reason, this set
is in one-to-one correspondence with the set of all permutations of size $n - 1$. A possible
bijection sends $w \in \SetCliff_\MapArithmetic{1}(n)$ to the permutation $\ToPermutation(w)
:= \sigma$ of size $n - 1$ where $\sigma$ is the permutation such that for any $i \in [n -
1]$, there are in $\sigma$ exactly $w(i)$ letters on the right of $i$ which are smaller than
$i$. The word $w$ is sometimes called the \Def{Lehmer code} of $\sigma$~\cite{Leh60}, up to
a slight variation. For instance, in $\SetCliff_\MapArithmetic{1}(7)$,
\begin{equation}
    \ToPermutation(002323) = 436512.
\end{equation}
The poset $\Par{\SetCliff_\MapArithmetic{1}(n), \Leq}$ is studied in~\cite{Den13}.

\subsubsection{$m$-increasing trees} \label{subsubsec:increasing_trees}
For any $m \geq 0$ and any $n \geq 1$, $\SetCliff_\MapArithmetic{m}(n)$ is the set of all
words $w$ of length $n - 1$ such that $w(i) \in \HanL{(i - 1) m}$ for all $i \in [n - 1]$.
For any $m \geq 0$ and $n \geq 1$,
\begin{equation}
    \# \SetCliff_\MapArithmetic{m}(n) = \prod_{i \in [n - 1]} \Par{1 + (i - 1)m}.
\end{equation}
This set is in one-to-one correspondence with the set of all \Def{$m$-increasing trees} with
$n - 1$ internal nodes, which are planar rooted trees where internal nodes are bijectively
labeled from $1$ to $n - 1$, have $m + 1$ children, and the sequence of the labels of the
nodes of any path starting from the root to the leaves is increasing. A possible
bijection~\cite{CG22} sends $w \in \SetCliff_\MapArithmetic{m}(n)$ to the $m$-increasing
tree $\ToTree_m(w) := \TreeT$ defined recursively as follows. If $w = \epsilon$, then
$\TreeT$ is the leaf. Otherwise, $w$ decomposes as $w = w' a$ where $w' \in
\SetCliff_\MapArithmetic{m}(n - 1)$ and $a \in \N$. In this case, $\TreeT$ is obtained by
grafting on the $a + 1$-st leaf of $\ToTree_m\Par{w'}$ an internal node labeled by $n$. For
instance, in $\SetCliff_\MapArithmetic{2}(8)$,
\begin{equation}
    \ToTree_2(0230228) =
    \scalebox{.8}{
    \begin{tikzpicture}[Centering,xscale=0.18,yscale=0.10]
        \node[Leaf](0)at(0.00,-8.80){};
        \node[Leaf](11)at(7.00,-4.40){};
        \node[Leaf](12)at(8.00,-13.20){};
        \node[Leaf](14)at(9.00,-13.20){};
        \node[Leaf](15)at(10.00,-13.20){};
        \node[Leaf](17)at(11.00,-13.20){};
        \node[Leaf](19)at(12.00,-13.20){};
        \node[Leaf](2)at(1.00,-8.80){};
        \node[Leaf](20)at(13.00,-13.20){};
        \node[Leaf](21)at(14.00,-8.80){};
        \node[Leaf](3)at(2.00,-17.60){};
        \node[Leaf](5)at(3.00,-17.60){};
        \node[Leaf](6)at(4.00,-17.60){};
        \node[Leaf](8)at(5.00,-13.20){};
        \node[Leaf](9)at(6.00,-13.20){};
        \node[Node](1)at(1.00,-4.40){$4$};
        \node[Node](10)at(7.00,0.00){$1$};
        \node[Node](13)at(9.00,-8.80){$7$};
        \node[Node](16)at(12.00,-4.40){$2$};
        \node[Node](18)at(12.00,-8.80){$3$};
        \node[Node](4)at(3.00,-13.20){$6$};
        \node[Node](7)at(5.00,-8.80){$5$};
        \draw[Edge](0)--(1);
        \draw[Edge](1)--(10);
        \draw[Edge](11)--(10);
        \draw[Edge](12)--(13);
        \draw[Edge](13)--(16);
        \draw[Edge](14)--(13);
        \draw[Edge](15)--(13);
        \draw[Edge](16)--(10);
        \draw[Edge](17)--(18);
        \draw[Edge](18)--(16);
        \draw[Edge](19)--(18);
        \draw[Edge](2)--(1);
        \draw[Edge](20)--(18);
        \draw[Edge](21)--(16);
        \draw[Edge](3)--(4);
        \draw[Edge](4)--(7);
        \draw[Edge](5)--(4);
        \draw[Edge](6)--(4);
        \draw[Edge](7)--(1);
        \draw[Edge](8)--(7);
        \draw[Edge](9)--(7);
        \node(r)at(7.00,4){};
        \draw[Edge](r)--(10);
    \end{tikzpicture}}.
\end{equation}
The posets $\Par{\SetCliff_{\MapArithmetic{m}}(n), \Leq}$ are studied in~\cite{CG20,CG22}.

\subsubsection{$c$-rectangular paths} \label{subsubsec:rectangular_paths}
For any $c \geq 0$ and any $n \geq 1$, $\SetHill_{\MapConstant{c}}(n)$ is the set of all
words $w$ of length $n - 1$ such that
\begin{math}
    w = 0^{\alpha_0} 1^{\alpha_1} 2^{\alpha_2} \dots c^{\alpha_c}
\end{math}
for a sequence $\Par{\alpha_0, \alpha_1, \alpha_2, \dots, \alpha_c}$ of nonnegative integers
such that $\alpha_0 + \alpha_1 + \alpha_2 + \dots + \alpha_c = n - 1$. It is straightforward
to show that
\begin{equation} \label{equ:number_rectangles}
    \# \SetHill_\MapConstant{c}(n) = \binom{n + c - 1}{c}.
\end{equation}
This set is in one-to-one correspondence with the set of all \Def{$c$-rectangular paths} of
size $n$ that are paths from $(0, 0)$ to $(n - 1, c)$ made of east steps $(1, 0)$ and north
steps $(0, 1)$. A possible bijection sends $w \in \SetHill_{\MapConstant{c}}(n)$ to the
$c$-rectangular path $\ToPath_c(w)$ having $\alpha_i$ east steps at ordinate~$i$. For
instance, in $\SetHill_{\MapConstant{4}}(8)$,
\begin{equation}
    \ToPath_4(1111244) =
    \scalebox{.8}{
    \begin{tikzpicture}[Centering,scale=.35]
        \draw[PathGrid](0,0)grid(7,4);
        \node[PathNode](0)at(0,0){};
        \node[PathNode](1)at(0,1){};
        \node[PathNode](2)at(1,1){};
        \node[PathNode](3)at(2,1){};
        \node[PathNode](4)at(3,1){};
        \node[PathNode](5)at(4,1){};
        \node[PathNode](6)at(4,2){};
        \node[PathNode](7)at(5,2){};
        \node[PathNode](8)at(5,3){};
        \node[PathNode](9)at(5,4){};
        \node[PathNode](10)at(6,4){};
        \node[PathNode](11)at(7,4){};
        \draw[PathStep](0)--(1);
        \draw[PathStep](1)--(2);
        \draw[PathStep](2)--(3);
        \draw[PathStep](3)--(4);
        \draw[PathStep](4)--(5);
        \draw[PathStep](5)--(6);
        \draw[PathStep](6)--(7);
        \draw[PathStep](7)--(8);
        \draw[PathStep](8)--(9);
        \draw[PathStep](9)--(10);
        \draw[PathStep](10)--(11);
    \end{tikzpicture}}.
\end{equation}
In the posets $\Par{\SetHill_{\MapConstant{c}}(n), \Leq}$, one has $u \Leq v$ if and only if
the $c$-rectangular path $\ToPath_c(v)$ is weakly above $\ToPath_c(u)$.

\subsubsection{$m$-Dyck paths} \label{subsubsec:dyck_paths}
For any $m \geq 0$ and any $n \geq 1$, $\SetHill_{\MapArithmetic{m}}(n)$ is the set of all
weakly increasing words $w$ of length $n - 1$ such that $w(i) \in \HanL{(i - 1) m}$ for all
$i \in [n - 1]$. It is shown in~\cite{CG20,CG22} that
\begin{math}
    \# \SetHill_\MapArithmetic{m}(n) = \FussCatalan_m(n - 1)
\end{math}
where
\begin{equation}
    \FussCatalan_m(n) := \frac{1}{mn + 1}\binom{mn + n}{n}
\end{equation}
is the $n$-th $m$-Fuss-Catalan number~\cite{DM47}. This set is in one-to-one correspondence
with the set of all \Def{$m$-Dyck paths} of size $n - 1$, that are paths from $(0, 0)$ to
$((m + 1) n, 0)$ staying above the $x$-axis and consisting in steps $(1, m)$ and $(1, -1)$.
A possible bijection sends $w \in \SetHill_{\MapArithmetic{m}}(n)$ to the $m$-Dyck path
$\ToDyck_m(w)$ such that for any $i \in [n - 1]$, the $i$-th step $(1, m)$ of $\ToDyck_m(w)$
has $w(i)$ steps $(1, -1)$ on its left. For instance, in $\SetHill_{\MapArithmetic{2}}(6)$,
\begin{equation}
    \ToDyck_2(02366) =
    \scalebox{.8}{
    \begin{tikzpicture}[Centering,scale=.3]
        \draw[PathGrid](0,0)grid(15,4);
        \node[PathNode](0)at(0,0){};
        \node[PathNode](1)at(1,2){};
        \node[PathNode](2)at(2,1){};
        \node[PathNode](3)at(3,0){};
        \node[PathNode](4)at(4,2){};
        \node[PathNode](5)at(5,1){};
        \node[PathNode](6)at(6,3){};
        \node[PathNode](7)at(7,2){};
        \node[PathNode](8)at(8,1){};
        \node[PathNode](9)at(9,0){};
        \node[PathNode](10)at(10,2){};
        \node[PathNode](11)at(11,4){};
        \node[PathNode](12)at(12,3){};
        \node[PathNode](13)at(13,2){};
        \node[PathNode](14)at(14,1){};
        \node[PathNode](15)at(15,0){};
        \draw[PathStep](0)--(1);
        \draw[PathStep](1)--(2);
        \draw[PathStep](2)--(3);
        \draw[PathStep](3)--(4);
        \draw[PathStep](4)--(5);
        \draw[PathStep](5)--(6);
        \draw[PathStep](6)--(7);
        \draw[PathStep](7)--(8);
        \draw[PathStep](8)--(9);
        \draw[PathStep](9)--(10);
        \draw[PathStep](10)--(11);
        \draw[PathStep](11)--(12);
        \draw[PathStep](12)--(13);
        \draw[PathStep](13)--(14);
        \draw[PathStep](14)--(15);
    \end{tikzpicture}}.
\end{equation}
The posets $\Par{\SetHill_{\MapArithmetic{1}}(n), \Leq}$ are the Stanley
lattices~\cite{Sta75,Knu06}. For $m \geq 2$, the posets
$\Par{\SetHill_{\MapArithmetic{m}}(n), \Leq}$ are generalizations of the previous
lattices~\cite{CG20,CG22}.

\subsection{Operads and interstice operads} \label{subsec:operads}
Let us provide some elementary definitions about nonsymmetric operads and interstice
operads.

\subsubsection{Graded spaces}
A \Def{graded space} is a set expressed as a direct sum
\begin{equation}
    \SpaceV := \bigoplus_{n \in \N} \SpaceV(n)
\end{equation}
such that all $\SpaceV(n)$, $n \in \N$, are spaces. Given $f \in \SpaceV$, if there is an $n
\in \N$ such that $f \in \SpaceV(n)$, $f$ is \Def{homogeneous}. In this case, this $n$ is
unique and is the \Def{rank} $|f|$ of $f$. If all $\SpaceV(n)$, $n \in \N$, have finite
dimensions, the \Def{Hilbert series} of $\SpaceV$ is the generating series $\sum_{n \in \N}
\dim \SpaceV(n) t^n$. If $\SpaceV$ and $\SpaceV'$ are two graded spaces, a linear map $\phi
: \SpaceV \to \SpaceV'$ is a \Def{graded space morphism} if for any $n \in \N$ and any $f
\in \SpaceV(n)$, $\phi(f) \in \SpaceV'(n)$. Besides, $\SpaceV'$ is a \Def{graded subspace}
of $\SpaceV$ if for any $n \in \N$, $\SpaceV'(n) \subseteq \SpaceV(n)$. Given a graded set
$S$, the linear span $\Span(S)$ of $S$ is the graded space defined as the direct sum of the
linear spans of each $S(n)$, $n \in \N$. By definition, the bases of $\Span(S)$ are indexed
by $S$. The \Def{elementary basis} (or \Def{$\BasisE$-basis} for short) of $\Span(S)$ is the
set $\Bra{\BasisE_x : x \in S}$ where each $\BasisE_x$, $x \in S$ is a formal symbol. We
shall consider several bases of a same graded space $\Span(S)$ defined from the
$\BasisE$-basis.

\subsubsection{Nonsymmetric operads}
A \Def{nonsymmetric operad in the category of spaces}, or a \Def{nonsymmetric operad} for
short, is a graded space $\Operad$ together with maps
\begin{equation}
    \Composition_i : \Operad(n) \otimes \Operad(m) \to \Operad(n + m - 1),
    \qquad 1 \leq i \leq n, \enspace 1 \leq m,
\end{equation}
called \Def{partial compositions}, and a distinguished element $\Unit \in \Operad(1)$, the
\Def{unit} of $\Operad$. These data have to satisfy, for any homogeneous elements $f_1$,
$f_2$,and $f_3$ of $\Operad$, the three relations
\begin{subequations}
\begin{equation} \label{equ:operad_axiom_1}
    \Par{f_1 \Composition_i f_2} \Composition_{i + j - 1} f_3
    = f_1 \Composition_i \Par{f_2 \Composition_j f_3},
    \qquad i \in \Han{\Brr{f_1}}, \enspace j \in \Han{\Brr{f_2}},
\end{equation}
\begin{equation} \label{equ:operad_axiom_2}
    \Par{f_1 \Composition_i f_2} \Composition_{j + \Brr{f_2} - 1} f_3
    = \Par{f_1 \Composition_j f_3} \Composition_i f_2,
    \qquad i, j \in \Han{\Brr{f_1}}, \enspace i < j,
\end{equation}
\begin{equation} \label{equ:operad_axiom_3}
    \Unit \Composition_1 f = f = f \Composition_i \Unit,
    \qquad i \in [|f|].
\end{equation}
\end{subequations}

We use in this work the definitions and conventions about nonsymmetric operads presented
in~\cite[Chapter 5]{Gir18}. Since this work deals only with nonsymmetric operads, we call
them simply \Def{operads}. Other usual references about operads are~\cite{Men15} for a
combinatorial point of view and~\cite{LV12} for an algebraic one. We shall use notions like
set-operads, free operads, ideals and quotients, presentations by generators and relations,
and Koszul duality.

\subsubsection{Interstice operads} \label{subsubsec:interstice_operads}
Given a set $A$, we see the set $A^*$ of the words on $A$ as a graded set such that the size
of $w \in A^*$ is $\Length(w) + 1$. Let us define on the graded set $\IntersticeOperad(A) :=
\Span\Par{A^*}$ the partial composition maps $\Composition_i$ defined linearly on the
$\BasisE$-basis of $\IntersticeOperad(A)$, for any $u, v \in A^*$ and $i \in[|u|]$, by
\begin{math}
    \BasisE_u \Composition_i \BasisE_v := \BasisE_{u \WhiteSquare_i v},
\end{math}
where
\begin{equation}
    u \WhiteSquare_i v := u(1, i - 1) \ v \ u(i, \Length(u)).
\end{equation}
For instance, in $\IntersticeOperad(\{\GenA, \GenB, \GenC\})$ we have
\begin{equation}
    \BasisE_{\ColB{\GenA \GenA \GenB \GenA \GenC \GenB}}
    \Composition_4
    \BasisE_{\ColA{\GenC \GenB \GenA \GenA}}
    = \BasisE_{\ColB{\GenA \GenA \GenB} \, \ColA{\GenC \GenB \GenA \GenA} \,
    \ColB{\GenA \GenC \GenB}}.
\end{equation}
It is straightforward to check that this structure is an operad admitting moreover
$\BasisE_\epsilon$ as unit. We call $\IntersticeOperad(A)$ the \Def{$A$-interstice operad}.

This operad is generated by the set $\GeneratingSet := \Bra{\BasisE_a : a \in A}$ of binary
elements. These generators are subjected exactly to the nontrivial relations
\begin{equation}
    \BasisE_b \Composition_1 \BasisE_a - \BasisE_a \Composition_2 \BasisE_b
\end{equation}
for all $a, b \in A$. In other terms, all elements of $A$ are binary generators, and these
elements are all associative with respect to all others. In particular, the algebras over
the operad $\IntersticeOperad(\{a, b\})$ are known under the name of \Def{duplexes of
vertices of cubes}~\cite[Section 6.3]{Pir03} (see also~\cite[Section 3]{AL07}
and~\cite[Definition 3.8]{ZGG20}).

\section{Operads of cliffs} \label{sec:cliff_operads}
By construction, $\IntersticeOperad(\N)$ is an operad on the linear span of the set all
words of nonnegative integers. Our aim is to build a substructure of $\IntersticeOperad(\N)$
on the linear span of $\SetCliff_\delta$ for the largest possible class of range maps
$\delta$. We propose a construction in the case where $\delta$ is unimodal. By defining
alternative bases of the obtained operads, we shall prove that these operads are set-operads
and provide some properties about their generators and their nontrivial relations.

\subsection{A quotient of an interstice operad} \label{subsec:quotient_interstice_operad}
We detail here the construction of the operads $\SpaceCliff_\delta$ when $\delta$ are
unimodal range maps.

\subsubsection{For weakly increasing range maps}
\label{subsubsec:operad_weakly_increasing_range_map}
First of all, if $\delta$ is a weakly increasing range map, it is immediate that the linear
span of the set $\Bra{\BasisE_w : w \in \SetCliff_\delta}$ forms a suboperad of
$\IntersticeOperad(\N)$. Let us denote by $\SpaceCliff_\delta$ this operad.

\subsubsection{Failure of direct quotients}
Observe that, given a (even unimodal) range map $\delta$, the graded subspace
$\SpaceV_\delta'$ of $\IntersticeOperad(\N)$ defined as the linear span of the set
$\Bra{\BasisE_w : w \in \N^* \setminus \SetCliff_\delta}$ is not always an operad ideal of
$\IntersticeOperad(\N)$. Indeed, for $\delta := 0110^\omega$, one has $\BasisE_{11} \in
\SpaceV_\delta'$ and $\BasisE_0 \in \IntersticeOperad(\N)$, but the element $\BasisE_{11}
\Composition_1 \BasisE_0 = \BasisE_{011}$ is not in~$\SpaceV_\delta'$. As a matter of fact,
it is possible to prove that $\SpaceV_\delta'$ is an operad ideal of $\IntersticeOperad(\N)$
if and only if $\delta$ is weakly increasing.

\subsubsection{For unimodal range maps}
The key of the construction is not to consider quotients of $\IntersticeOperad(\N)$ but
quotients of some suboperads of $\IntersticeOperad(\N)$ instead.

For this, given any range map $\delta$, we denote by $\bar{\delta}$ the range map defined by
$\bar{\delta}(i) := \max \Bra{\delta(1), \dots, \delta(i)}$ for any $i \geq 1$. For
instance, if $\delta = 10032242^\omega$, then $\bar{\delta} = 11133344^\omega$. Remark that
when $\delta$ is weakly increasing, $\bar{\delta} = \delta$. By construction, the range map
$\bar{\delta}$ is weakly increasing. For this reason, $\SpaceCliff_{\bar{\delta}}$ is a
well-defined operad. Observe moreover that $\SetCliff_\delta$ is a subset of
$\SetCliff_{\bar{\delta}}$. Let $\SpaceV_\delta$ be the graded subspace of
$\SpaceCliff_{\bar{\delta}}$ defined as the linear span of the set
\begin{math}
    \Bra{\BasisE_w : w \in \SetCliff_{\bar{\delta}} \setminus \SetCliff_\delta}.
\end{math}

\begin{Proposition} \label{prop:unimodality_operad_ideal_cliffs}
    For any range map $\delta$, the space $\SpaceV_\delta$ is an operad ideal of the operad
    $\SpaceCliff_{\bar{\delta}}$ if and only if $\delta$ is unimodal.
\end{Proposition}
\begin{proof}
    Assume first that $\delta$ is not unimodal. Thus, there are indices $1 \leq \alpha_1 <
    \alpha_2 < \alpha_3$ such that $\delta\Par{\alpha_1} > \delta\Par{\alpha_2} <
    \delta\Par{\alpha_3}$. Let $u := 0^{\alpha_2 - 1} \Par{\delta\Par{\alpha_2} + 1}$. By
    construction, $u \notin \SetCliff_\delta$ and, since $\delta\Par{\alpha_1} \geq
    \delta\Par{\alpha_2} + 1$, $u \in \SetCliff_{\bar{\delta}}$. Therefore, $\BasisE_u \in
    \SpaceV_\delta$. Let also $v := 0^{\alpha_3 - \alpha_2}$ so that $\BasisE_v \in
    \SpaceCliff_{\bar{\delta}}$. Now, one has $\BasisE_u \Composition_1 \BasisE_v =
    \BasisE_w$ with $w := 0^{\alpha_3 - 1} \Par{\delta\Par{\alpha_2} + 1}$. Since
    $\delta\Par{\alpha_3} \geq \delta\Par{\alpha_2} + 1$, we have $w \in \SetCliff_\delta$,
    so that $\BasisE_w \notin \SpaceV_\delta$. This shows that $\SpaceV_\delta$ is not an
    operad ideal of $\SpaceCliff_{\bar{\delta}}$.

    Conversely, assume that $\delta$ is unimodal. Let $u, v \in \SetCliff_{\bar{\delta}}$
    and $i \in [|u|]$ with $u \notin \SetCliff_\delta$ or $v \notin \SetCliff_\delta$. Let
    us set $\BasisE_w := \BasisE_u \Composition_i \BasisE_v$. From these assumptions, we
    obtain that there is an index $j \geq 1$ and a letter $a$ at position $j$ of $u$ or of
    $v$ such that $a > \delta(j)$. Since $u \in \SetCliff_{\bar{\delta}}$ and $v \in
    \SetCliff_{\bar{\delta}}$, there is an index $j' < j$ such that $\delta\Par{j'} >
    \delta(j)$ and $a \leq \delta\Par{j'}$. Moreover, since $\delta$ is unimodal, for all
    $j'' \geq j$, $\delta(j) \geq \delta\Par{j''}$. Now, due to the definition of the
    partial composition of $\SetCliff_{\bar{\delta}}$, the letter $a$ appears in $w$ at a
    certain position $j + k$ for a $k \geq 0$. Since
    \begin{math}
        w(j + k) = a > \delta(j) \geq \delta(j + k),
    \end{math}
    we have that $w \notin \SetCliff_\delta$ and thus $\BasisE_w \in \SpaceV_\delta$.
    Therefore, and since $\BasisE_\epsilon \notin \SpaceV_\delta$, $\SpaceV_\delta$ is an
    operad ideal of~$\SpaceCliff_{\bar{\delta}}$.
\end{proof}

As a consequence of Proposition~\ref{prop:unimodality_operad_ideal_cliffs}, one has the
following result.

\begin{Theorem} \label{thm:cliff_operad_unimodal_delta}
    For any unimodal range map $\delta$, the space $\SpaceCliff_{\bar{\delta}} /
    \SpaceV_\delta$ is an operad.
\end{Theorem}

For any unimodal range map $\delta$, we set $\SpaceCliff_\delta :=
\SpaceCliff_{\bar{\delta}} / \SpaceV_\delta$. Since when $\delta$ is weakly increasing,
$\SpaceV_\delta$ is the null space, this definition is consistent with the previous
definition of $\SpaceCliff_\delta$ given in
Section~\ref{subsubsec:operad_weakly_increasing_range_map}. By construction, the partial
composition of $\SpaceCliff_\delta$ satisfies, for any $u, v \in \SetCliff_\delta$ and $i
\in [|u|]$,
\begin{equation} \label{equ:composition_e_basis}
    \BasisE_u \Composition_i \BasisE_v
    = \CharacteristicCliff_\delta\Par{u \WhiteSquare_i v} \BasisE_{u \WhiteSquare_i v},
\end{equation}
where $\CharacteristicCliff_\delta : \N^* \to \K$ is the map defined for any $w \in \N^*$ by
$\CharacteristicCliff_\delta(w) := 1$ if $w \in \SetCliff_\delta$ and by
$\CharacteristicCliff_\delta(w) := 0$ otherwise. For instance, in
$\SpaceCliff_{1232^\omega}$ we have
\begin{subequations}
\begin{equation}
    \BasisE_{\ColB{002}} \Composition_3 \BasisE_{\ColA{10}}
    = \BasisE_{\ColB{00} \ColA{10} \ColB{2}},
\end{equation}
\begin{equation}
    \BasisE_{\ColB{002}} \Composition_3 \BasisE_{\ColA{1311}} = 0.
\end{equation}
\end{subequations}

\subsection{Fundamental and homogeneous bases}
The aim of this section is to introduce two alternative bases of $\SpaceCliff_\delta$ in
order to prove that this operad is a set-operad. The first of these is the fundamental
basis, which is defined from the elementary basis by using the poset structure of
$\delta$-cliffs. The second is the homogeneous basis which is defined from the fundamental
basis, again by using the same poset structure but in a different way. The partial
composition of $\SpaceCliff_\delta$ is expressed here over these two alternative bases.

\subsubsection{Fundamental basis}
Let $\delta$ be a unimodal range map. For any $u \in \SetCliff_\delta$, let
\begin{equation}
    \BasisF_w := \sum_{\substack{w' \in \SetCliff_\delta \\ w \Leq w'}}
    \MobiusFunction_\Leq\Par{w, w'} \BasisE_{w'},
\end{equation}
where $\MobiusFunction_\Leq$ is the Möbius function of the $\delta$-cliff posets. Since this
poset a Cartesian product of total orders, for any $w, w' \in \SetCliff_\delta$ such that $w
\Leq w'$,
\begin{math}
    \MobiusFunction_\Leq\Par{w, w'}
    =
    \prod_{i \in [\Length(w)]} \MobiusFunction\Par{w(i), w'(i)}
\end{math}
where, for any $a \leq a' \in \N$,
\begin{equation}
    \MobiusFunction\Par{a, a'}
    :=
    \begin{cases}
        1 & \mbox{if } a' = a, \\
        -1 & \mbox{if } a' = a + 1, \\
        0 & \mbox{otherwise}.
    \end{cases}
\end{equation}
For instance, in $\SpaceCliff_{224^\omega}$ we have
\begin{equation}
    \BasisF_{1221}
        = \BasisE_{1221}
        - \BasisE_{1222} - \BasisE_{1231} - \BasisE_{2221}
        + \BasisE_{1232} + \BasisE_{2222} + \BasisE_{2231} 
        - \BasisE_{2232}.
\end{equation}
By Möbius inversion and triangularity, for any $u \in \SetCliff_\delta$,
\begin{equation}
    \BasisE_w = \sum_{\substack{w' \in \SetCliff_\delta \\ w \Leq w'}} \BasisF_{w'},
\end{equation}
so that the set $\Bra{\BasisF_w : w \in \SetCliff_\delta}$ is a basis of
$\SpaceCliff_\delta$, called \Def{fundamental basis} (or \Def{$\BasisF$-basis} for short).

Given $w \in \SetCliff_\delta$ and a map $f : [\Length(w)] \to \N \setminus \{0\}$, let
$\Modify_\delta^f(w)$ be the word on $\N$ satisfying, for any $j \in [\Length(w)]$,
\begin{equation}
    \Par{\Modify_\delta^f(w)}(j) =
    \begin{cases}
        \delta(f(j)) & \mbox{if } w(j) = \delta(j), \\
        w(j) & \mbox{otherwise}.
    \end{cases}
\end{equation}
For any $u, v \in \SetCliff_\delta$ and $i \in [|u|]$, let
\begin{equation}
    u \BlackSquare_i v := \Modify_\delta^f(u) \WhiteSquare_i \Modify_\delta^g(v)
\end{equation}
where $f(j) := j$ if $j \leq i - 1$ and $f(j) := j + \Length(v)$ otherwise, and $g(j) := j +
i - 1$. For instance for $\delta := 11321^\omega$, one has
\begin{subequations}
\begin{equation} \label{equ:example_black_composition}
    \ColB{{\bf 1} 02 {\bf 2}} \BlackSquare_3 \ColA{{\bf 1} 01}
    = \ColB{{\bf 1} 0} \ \ColA{{\bf 3} 01} \ \ColB{2 {\bf 1}},
\end{equation}
\begin{equation}
    \ColB{{\bf 1} 02 {\bf 2}} \BlackSquare_4 \ColA{00 {\bf 3}}
    = \ColB{{\bf 1} 0 2} \ \ColA{00 {\bf 1}} \ \ColB{{\bf 1}}.
\end{equation}
\end{subequations}
Observe from~\eqref{equ:example_black_composition} that even if $1022$ and $101$ are two
$\delta$-cliffs, $1030121$ is not a $\delta$-cliff.

\begin{Lemma} \label{lem:good_definition_black_composition}
    Let $\delta$ be a range map, $u, v \in \SetCliff_\delta$, and $i \in [|u|]$. If $u
    \WhiteSquare_i v$ is a $\delta$-cliff, then $u \BlackSquare_i v$ also is.
\end{Lemma}
\begin{proof}
    Let $w := u \WhiteSquare_i v$ and assume that $w \in \SetCliff_\delta$. Thus, for all $j
    \in [\Length(w)]$, $w(j) \leq \delta(j)$. By setting $w' := u \BlackSquare_i v$, by
    definition of the operation $\BlackSquare_i$ and the previous hypothesis, $w'(j) \leq
    \max \{w(j), \delta(j)\}$. Therefore, $w' \in \SetCliff_\delta$.
\end{proof}

\begin{Lemma} \label{lem:interval_white_and_black_squares}
    Let $\delta$ be a range map, $u, v \in \SetCliff_\delta$, and $i \in [|u|]$ such that $u
    \WhiteSquare_i v$ is a $\delta$-cliff. We have $w \in \Han{u \WhiteSquare_i v, u
    \BlackSquare_i v}_\Leq$ if and only if the following three assertions hold.
    \begin{enumerate}[label={\it (\roman*)}]
        \item For any $j \in [i - 1]$, $w(j) = u(j)$.
        \item For any $j \in [i, \Length(v) + i - 1]$, if $v(j - i + 1) = \delta(j - i + 1)$
        then $w(j) \in [v(j - i + 1), \delta(j)]$, and $w(j) = v(j - i + 1)$ otherwise.
        \item For any $j \in [\Length(v) + i, \Length(w)]$, if $u(j - \Length(v)) = \delta(j
        - \Length(v))$ then $w(j) \in [u(j - \Length(v)), \delta(j)]$, and $w(j) = u(j -
        \Length(v))$ otherwise.
    \end{enumerate}
\end{Lemma}
\begin{proof}
    This is a direct consequence of the definitions of the operations $\WhiteSquare_i$ and
    $\BlackSquare_i$.
\end{proof}

\begin{Proposition} \label{prop:composition_f_basis}
    For any unimodal range map $\delta$, the partial composition map $\Composition_i$ on
    $\SpaceCliff_\delta$ satisfies, for any $u, v \in \SetCliff_\delta$ and $i \in [|u|]$,
    \begin{equation} \label{equ:composition_f_basis}
        \BasisF_u \Composition_i \BasisF_v
        = \CharacteristicCliff_\delta\Par{u \WhiteSquare_i v}
        \sum_{w \in \Han{u \WhiteSquare_i v, u \BlackSquare_i v}_\Leq}
        \BasisF_w.
    \end{equation}
\end{Proposition}
\begin{proof}
    Let us denote by $\Composition'_i$ the operation on $\SpaceCliff_\delta$ defined
    in~\eqref{equ:composition_f_basis} and let us show that this operation is the same as
    $\Composition_i$.  First of all, by Lemma~\ref{lem:good_definition_black_composition},
    for any $u, v \in \SetCliff_\delta$, if $\CharacteristicCliff_\delta\Par{u
    \WhiteSquare_i v} = 1$, then $\Han{u \WhiteSquare_i v, u \BlackSquare_i v}_\Leq$ is an
    interval of a $\delta$-cliff poset. For this reason,
    \eqref{equ:composition_f_basis} is well-defined. By Möbius inversion,
    \begin{equation} \label{equ:composition_f_basis_1}
        \BasisE_u \Composition'_i \BasisE_v
        = \sum_{\substack{
            u', v' \in \SetCliff_\delta \\
            u \Leq u' \\
            v \Leq v'
        }}
        \BasisF_{u'} \Composition'_i \BasisF_{v'}
        =
        \sum_{\substack{
            u', v' \in \SetCliff_\delta \\
            u \Leq u' \\
            v \Leq v'
        }}
        \CharacteristicCliff_\delta\Par{u' \WhiteSquare_i v'}
        \sum_{w \in \Han{u' \WhiteSquare_i v', u' \BlackSquare_i v'}_\Leq}
        \BasisF_w.
    \end{equation}
    Let $z := u \WhiteSquare_i v$.

    Assume first that $z \notin \SetCliff_\delta$. Thus, there is $j \in [\Length(z)]$ such
    that $z(j) > \delta(j)$. For all $u', v' \in \SetCliff_\delta$ such that $u \Leq u'$ and
    $v \Leq v'$, by setting $z' := u' \WhiteSquare_i v'$, one has $z \Leq z'$. Therefore, we
    have $\delta(j) < z(j) \leq z'(j)$, showing that $z' \notin \SetCliff_\delta$.
    By~\eqref{equ:composition_f_basis_1}, $\BasisE_u \Composition'_i \BasisE_v = 0$ in this
    case.

    Assume now that $z \in \SetCliff_\delta$. First of all, we immediately have that if a
    $\BasisF_w$ appears in~\eqref{equ:composition_f_basis_1} with $w \in \SetCliff_\delta$,
    then $u \WhiteSquare_i v \Leq w$. Conversely, let $w \in \SetCliff_\delta$ such that $u
    \WhiteSquare_i v \Leq w$. By Lemma~\ref{lem:interval_white_and_black_squares}, there
    exists a unique pair $\Par{u', v'} \in \SetCliff_\delta^2$ such that $u \Leq u'$, $v
    \Leq v'$, and $w \in \Han{u' \WhiteSquare_i v', u' \BlackSquare_i v'}_\Leq$.

    From all this, we deduce by~\eqref{equ:composition_e_basis} that
    \begin{equation}
        \BasisE_u \Composition'_i \BasisE_v
        = \CharacteristicCliff_\delta\Par{u \WhiteSquare_i v}
        \sum_{\substack{
            w \in \SetCliff_\delta \\
            u \WhiteSquare_i v \Leq w
        }}
        \BasisF_w
        =
        \CharacteristicCliff_\delta\Par{u \WhiteSquare_i v}
        \BasisE_{u \WhiteSquare_i v}
        =
        \BasisE_u \Composition_i \BasisE_v.
    \end{equation}
    Therefore, the operations $\Composition'_i$ and $\Composition_i$ are the same,
    establishing~\eqref{equ:composition_f_basis}.
\end{proof}

For instance, in $\SpaceCliff_{123454^\omega}$ we have
\begin{subequations}
\begin{equation}
    \BasisF_{\ColB{{\bf 1} 0}} \Composition_2 \BasisF_{\ColA{0 {\bf 2} 1}} =
    \BasisF_{\ColB{1} \ColA{021} \ColB{0}} + \BasisF_{\ColB{1} \ColA{031} \ColB{0}},
\end{equation}
\begin{equation}
    \BasisF_{\ColB{01 {\bf 3}}} \Composition_2 \BasisF_{\ColA{{\bf 1} 0 {\bf 3}}}
    = \BasisF_{\ColB{0} \ColA{103} \ColB{13}}
    + \BasisF_{\ColB{0} \ColA{103} \ColB{14}}
    + \BasisF_{\ColB{0} \ColA{104} \ColB{13}}
    + \BasisF_{\ColB{0} \ColA{104} \ColB{14}}
    + \BasisF_{\ColB{0} \ColA{203} \ColB{13}}
    + \BasisF_{\ColB{0} \ColA{203} \ColB{14}}
    + \BasisF_{\ColB{0} \ColA{204} \ColB{13}}
    + \BasisF_{\ColB{0} \ColA{204} \ColB{14}}.
\end{equation}
\end{subequations}

\subsubsection{Homogeneous basis}
For any $w \in \SetCliff_\delta$, let
\begin{equation}
    \BasisH_w :=
    \sum_{\substack{w' \in \SetCliff_\delta \\
        w' \Leq w
    }}
    \BasisF_{w'}.
\end{equation}
For instance, in $\SpaceCliff_{3221^\omega}$ we have
\begin{equation}
    \BasisH_{2101}
    = \BasisF_{0000} + \BasisF_{0001} + \BasisF_{0101} + \BasisF_{1001} + \BasisF_{1100}
    + \BasisF_{1101} + \BasisF_{2000} + \BasisF_{2001} + \BasisF_{2100} + \BasisF_{2101}.
\end{equation}
By triangularity, the set $\Bra{\BasisH_w : w \in \SetCliff_\delta}$ is a basis of
$\SpaceCliff_\delta$, called \Def{homogeneous basis} (or \Def{$\BasisH$-basis} for short).

\begin{Proposition} \label{prop:composition_h_basis}
    For any unimodal range map $\delta$, the partial composition map $\Composition_i$ on
    $\SpaceCliff_\delta$ satisfies, for any $u, v \in \SetCliff_\delta$ and $i \in [|u|]$,
    \begin{equation}
        \BasisH_u \Composition_i \BasisH_v
        = \BasisH_{\Reduction_\delta\Par{u \BlackSquare_i v}}.
    \end{equation}
\end{Proposition}
\begin{proof}
    By Proposition~\ref{prop:composition_f_basis},
    \begin{equation} \label{equ:composition_h_basis_1}
        \BasisH_u \Composition_i \BasisH_v
        = \sum_{\substack{
            u', v' \in \SetCliff_\delta \\
            u' \Leq u \\
            v' \Leq v
        }}
        \BasisF_{u'} \Composition_i \BasisF_{v'}
        =
        \sum_{\substack{
            u', v' \in \SetCliff_\delta \\
            u' \Leq u \\
            v' \Leq v
        }}
        \CharacteristicCliff_\delta\Par{u' \WhiteSquare_i v'}
        \sum_{w \in \Han{u' \WhiteSquare_i v', u' \BlackSquare_i v'}_\Leq}
        \BasisF_w.
    \end{equation}
    First of all, we immediately have that if a $\BasisF_w$ appears
    in~\eqref{equ:composition_h_basis_1} with $w \in \SetCliff_\delta$, then $w \Leq
    \Reduction_\delta\Par{u \BlackSquare_i v}$.  Conversely, let $w \in \SetCliff_\delta$
    such that $w \Leq \Reduction_\delta\Par{u \BlackSquare_i v}$.  By
    Lemma~\ref{lem:interval_white_and_black_squares}, there exists a unique pair $\Par{u',
    v'} \in \SetCliff_\delta$ such that $u' \Leq u$, $v' \Leq v$, and $w \in \Han{u'
    \WhiteSquare_i v', u' \BlackSquare_i v'}_\Leq$. For this reason,
    from~\eqref{equ:composition_h_basis_1}, we obtain
    \begin{equation}
        \BasisH_u \Composition_i \BasisH_v
        =
        \sum_{\substack{
            w \in \SetCliff_\delta \\
            w \Leq \Reduction_\delta\Par{u \BlackSquare_i v}
        }}
        \BasisF_w
        =
        \BasisH_{\Reduction_\delta\Par{u \BlackSquare_i v}},
    \end{equation}
    showing the statement of the proposition.
\end{proof}

For instance, in $\SpaceCliff_{22342^\omega}$ we have
\begin{subequations}
\begin{equation}
    \BasisH_{\ColB{01}} \Composition_3 \BasisH_{\ColA{{\bf 2} {\bf 2} 1}}
    = \BasisH_{\ColB{01} \ColA{341}},
\end{equation}
\begin{equation}
    \BasisH_{\ColB{{\bf 2} 0 {\bf 3} 3}} \Composition_3 \BasisH_{\ColA{1 {\bf 2}}}
    = \BasisH_{\ColB{20} \ColA{14} \ColB{22}}.
\end{equation}
\end{subequations}

Recall that a basis of an operad is a \Def{set-operad basis} if any partial composition of
two basis elements is a basis element. An operad having a set-operad basis is a
\Def{set-operad}. As a consequence of Proposition~\ref{prop:composition_h_basis}, one has
the following result.

\begin{Theorem} \label{thm:set_operad}
    For any unimodal range map $\delta$, the operad $\SpaceCliff_\delta$ is a set-operad
    and its $\BasisH$-basis is a set-operad basis.
\end{Theorem}

\subsection{Generators and relations}
We provide here some results about minimal generating families of $\SpaceCliff_\delta$, the
existence of finite such families for these operads, and the nonexistence of finite families
of nontrivial relations for these operads. To explore the structure of the operads
$\SpaceCliff_\delta$, we work through the $\BasisE$-basis since the partial composition
expressed over this basis (see Equation~\eqref{equ:composition_e_basis}) is very elementary.

\subsubsection{Prime cliffs}
A nonempty $\delta$-cliff $w$ is \Def{$\delta$-prime} if the relation $w = u \WhiteSquare_i
v$ with $u, v \in \SetCliff_\delta$ and $i \in [|u|]$ implies $(u, v) \in \{(w, \epsilon),
(\epsilon, w)\}$. For instance, for $\delta := 1223321^\omega$, the $\delta$-cliffs $10033$
and $121332$ are $\delta$-prime, while $11222 = 122 \WhiteSquare_2 12$ is not. We denote by
$\SetPrime_\delta$ the graded subset of $\SetCliff_\delta$ of all $\delta$-prime
$\delta$-cliffs.

If $\delta$ is a nonconstant range map, let us denote by $\ConstantLength(\delta)$ the
smallest index $k \geq 1$ such that $\delta(k) \ne \delta(k + 1)$. For instance,
$\ConstantLength(\MapArithmetic{2}) = 1$ and $\ConstantLength(22{ \bf 2} 4445^\omega) = 3$.
If $\delta$ is $1$-dominated, we denote by $\DominationLength(\delta)$ the smallest index $k
\geq 1$ such that for all $k' \geq k$, $\delta(1) \geq \delta\Par{k'}$. For instance,
$\DominationLength\Par{\MapConstant{3}} = 1$ and $\DominationLength(22334 {\bf 1} 10^\omega)
= 6$.

\begin{Proposition} \label{prop:properties_prime_sets}
    Let $\delta$ be a unimodal range map.
    \begin{enumerate}[label={\it (\roman*)}]
        \item \label{item:properties_prime_sets_1}
        If $\delta$ is weakly decreasing, then $\SetPrime_\delta = \HanL{\delta(1)}$.
        \item \label{item:properties_prime_sets_2}
        Otherwise, if $\delta$ is $1$-dominated, then for any $w \in \SetPrime_\delta$,
        $\Length(w) \leq \DominationLength(\delta) - 1$.
        \item \label{item:properties_prime_sets_3}
        Otherwise, $\delta$ is not $1$-dominated and for any $k \geq \ConstantLength(\delta)
        + 1$, $\delta(1, k) \in \SetPrime_\delta$.
    \end{enumerate}
\end{Proposition}
\begin{proof}
    Assume that $\delta$ is weakly decreasing. It is immediate that any $\delta$-cliff
    having $1$ as length is $\delta$-prime. Moreover, let $w \in \SetCliff_\delta$ with
    $\Length(w) \geq 2$. We have $w = a w'$ with $a \in \HanL{\delta(1)}$ and $w'$ is a word
    of length $\Length(w) - 1$. Therefore, we have $w = a \WhiteSquare_2 w'$, and since
    $\delta$ is weakly decreasing, $w' \in \SetCliff_\delta$. This shows that $w \notin
    \SetPrime_\delta$ and implies~\ref{item:properties_prime_sets_1}.

    Assume that $\delta$ is not weakly decreasing and is $1$-dominated. Therefore
    $\DominationLength(\delta) \geq 2$. Let $w \in \SetCliff_\delta$ such that $\Length(w)
    \geq \DominationLength(\delta)$. One has
    \begin{equation}
        w =
        w(1, \DominationLength(\delta) - 1)
        \WhiteSquare_{\DominationLength(\delta)}
        w(\DominationLength(\delta)).
    \end{equation}
    By construction, $w(\DominationLength(\delta)) \leq \delta(\DominationLength(\delta))
    \leq \delta(1)$. For this reason, $w(\DominationLength(\delta))$ is a $\delta$-cliff.
    Moreover, since $w(1, \DominationLength(\delta) - 1)$ is a prefix of a $\delta$-cliff,
    this last word is also a $\delta$-cliff. This shows that $w$ is not $\delta$-prime and
    implies~\ref{item:properties_prime_sets_2}.

    Finally, assume that $\delta$ is not $1$-dominated and let $k \geq
    \ConstantLength(\delta) + 1$. Assume that there are $u, v \in \SetCliff_\delta$ and $i
    \in [|u|]$ such that $\delta(1, k) = u \WhiteSquare_i v$. The fact that there is an
    index $j \geq \ConstantLength(\delta) + 1$ such that $\delta(j) > \delta(1)$ implies
    that $u = \epsilon$ of $v = \epsilon$. Indeed, if $u \ne \epsilon$ and $v \ne \epsilon$,
    there would be an index $j' < j$ such that $u\Par{j'} = \delta(j)$ or $v\Par{j'} =
    \delta(j)$, contradicting the fact that $u$ and $v$ are $\delta$-cliffs. This
    implies~\ref{item:properties_prime_sets_3}.
\end{proof}

\subsubsection{Minimal generating sets}
Let us denote by $\GeneratingSet_\delta$ the family
$\Bra{\BasisE_w : w \in \SetPrime_\delta}$ of elements of $\SpaceCliff_\delta$.
For instance,
\begin{subequations}
\begin{equation} \label{equ:example_generating_family_1}
    \GeneratingSet_{21^\omega} = \Bra{\BasisE_0, \BasisE_1, \BasisE_2},
\end{equation}
\begin{equation} \label{equ:example_generating_family_2}
    \GeneratingSet_{1221^\omega} =
    \Bra{\BasisE_{0}, \BasisE_{1}, \BasisE_{02}, \BasisE_{12}, \BasisE_{022},
    \BasisE_{122}},
\end{equation}
\begin{multline} \label{equ:example_generating_family_3}
    \GeneratingSet_{11320^\omega} =
    \left\{\BasisE_{0}, \BasisE_{1}, \BasisE_{002}, \BasisE_{003}, \BasisE_{012},
    \BasisE_{013}, \BasisE_{102}, \BasisE_{103}, \BasisE_{112}, \BasisE_{113},
    \BasisE_{0022}, \BasisE_{0032},
    \right. \\ \left.
    \BasisE_{0122}, \BasisE_{0132}, \BasisE_{1022}, \BasisE_{1032}, \BasisE_{1122},
    \BasisE_{1132}\right\},
\end{multline}
\begin{equation} \label{equ:example_generating_family_4}
    \GeneratingSet_{112^\omega} =
    \Bra{\BasisE_{0}, \BasisE_{1}, \BasisE_{002}, \BasisE_{012}, \BasisE_{102},
    \BasisE_{112}, \BasisE_{0022}, \BasisE_{0122}, \BasisE_{1022}, \BasisE_{1122},
    \BasisE_{00222}, \BasisE_{01222}, \BasisE_{10222}, \BasisE_{11222}, \dots},
\end{equation}
\end{subequations}
By Proposition~\ref{prop:properties_prime_sets}, the
families~\eqref{equ:example_generating_family_1}, \eqref{equ:example_generating_family_2},
and ~\eqref{equ:example_generating_family_3} are finite, and the
family~\eqref{equ:example_generating_family_4} is infinite.

\begin{Proposition} \label{prop:minimal_generating_set_cliffs}
    For any unimodal range map $\delta$, $\GeneratingSet_\delta$ is a minimal generating set
    of the operad $\SpaceCliff_\delta$.
\end{Proposition}
\begin{proof}
    Let us first prove that $\GeneratingSet_\delta$ is a generating set of
    $\SpaceCliff_\delta$. By induction of the arity, it appears that any $w \in
    \SetCliff_\delta$ writes as an expression involving only $\delta$-prime cliffs and
    operations~$\WhiteSquare_i$. Therefore, due to the partial composition of
    $\SpaceCliff_\delta$ over the $\BasisE$-basis (see~\eqref{equ:composition_e_basis}), any
    $\BasisE_w$ writes as an expression involving only elements of $\GeneratingSet_\delta$
    and operations~$\Composition_i$. This shows the first assertion.

    Finally, the minimality of $\GeneratingSet_\delta$ follows from the fact each $\BasisE_w
    \in \GeneratingSet_\delta$ cannot be obtained as partial compositions of elements of
    $\GeneratingSet_\delta \setminus \Bra{\BasisE_w}$.
\end{proof}

\begin{Theorem} \label{thm:finiteness_generating_set_cliff_operads}
    For any unimodal range map $\delta$, the generating set $\GeneratingSet_\delta$ of the
    operad $\SpaceCliff_\delta$ is finite if and only if $\delta$ is $1$-dominated.
\end{Theorem}
\begin{proof}
    This is a consequence of Proposition~\ref{prop:minimal_generating_set_cliffs} and
    Points~\ref{item:properties_prime_sets_2} and~\ref{item:properties_prime_sets_3} of
    Proposition~\ref{prop:properties_prime_sets}.
\end{proof}

\subsubsection{Nontrivial relations}
Let us denote by $\RelationSpace_\delta$ the space of the nontrivial relations of
$\SpaceCliff_\delta$. We have here only a sufficient condition for the fact that
$\RelationSpace_\delta$ is not finitely generated.

\begin{Proposition} \label{prop:not_finiteness_relation_space_cliff_operads}
    For any unimodal range map $\delta$, if $\delta$ is not $1$-dominated, then the space
    $\RelationSpace_\delta$ is not finitely generated.
\end{Proposition}
\begin{proof}
    When $\delta$ is not $1$-dominated, by Point~\ref{item:properties_prime_sets_3} of
    Proposition~\ref{prop:properties_prime_sets}, for any $k \geq \ConstantLength(\delta) +
    1$, $\delta(1, k) \in \SetPrime_\delta$ so that $\BasisE_{\delta(1, k)} \in
    \GeneratingSet_\delta$.  Moreover, one has $\delta(1, k) \WhiteSquare_1 0 = 0 \delta(1,
    k) = 0 \WhiteSquare_2 \delta(1, k)$ so that
    \begin{equation}
        \BasisE_{\delta(1, k)} \Composition_1 \BasisE_0
        = \CharacteristicCliff_\delta(0 \delta(1, k)) \BasisE_{0 \delta(1, k)}
        = \BasisE_0 \Composition_2 \BasisE_{\delta(1, k)}.
    \end{equation}
    Since both $\BasisE_{\delta(1, k)}$ and $\BasisE_0$ belong to $\GeneratingSet_\delta$,
    the element $\BasisE_{\delta(1, k)} \Composition_1 \BasisE_0 - \BasisE_0 \Composition_2
    \BasisE_{\delta(1, k)}$ is an element of $\RelationSpace_\delta$. Since one has such a
    generating relation for each $k \geq \ConstantLength(\delta) + 1$, this implies the
    statement of the proposition.
\end{proof}

\section{Quotient operads} \label{sec:quotient_operads}
As exposed in Section~\ref{subsec:cliffs_and_objects}, some subsets $\SubFamilly$ of
$\delta$-cliffs satisfying some conditions are in one-to-one correspondence with other
graded sets (in particular, $c$-rectangular paths and $m$-Dyck paths). We describe here a
generic way to construct quotient operads of $\SpaceCliff_\delta$ whose bases are indexed
by~$\SubFamilly$.

\subsection{General construction}
We construct now a quotient of $\SpaceCliff_\delta$ by identifying some of its elements over
the $\BasisF$-basis with zero. In order to obtain a quotient operad, $\SubFamilly$ needs to
satisfy a condition which is stated now.

\subsubsection{Closure by subword reduction}
Let $\delta$ be a unimodal range map and $\SubFamilly$ be a nonempty graded subset of
$\SetCliff_\delta$. The graded set $\SubFamilly$ is \Def{closed by subword reduction} if for
any $w \in \SubFamilly$, for all subwords $w'$ of $w$ (that are sequences of not necessarily
contiguous letters of $w$), $\Reduction_\delta\Par{w'} \in \SubFamilly$. Remark that the
fact that $\SubFamilly$ is nonempty implies in this case that $\epsilon \in \SubFamilly$.

\subsubsection{Quotient operad}
Let the quotient space $\SpaceCliff_\SubFamilly := \SpaceCliff_\delta / \SpaceV_\SubFamilly$
where $\SpaceV_\SubFamilly$ is the linear span of the set $\Bra{\BasisF_w : w \in
\SetCliff_\delta \setminus \SubFamilly}$. This set is the \Def{fundamental basis} (or
\Def{$\BasisF$-basis} for short) of $\SpaceCliff_\SubFamilly$.

\begin{Proposition} \label{prop:quotient_cliff_operads}
    Let $\delta$ be a unimodal range map $\delta$ and $\SubFamilly$ be a nonempty graded
    subset of $\SetCliff_\delta$. If $\SubFamilly$ is closed by subword reduction, then
    $\SpaceV_\SubFamilly$ is an operad ideal of $\SpaceCliff_\delta$. Therefore, in this
    case, $\SpaceCliff_\SubFamilly$ is a quotient operad of $\SpaceCliff_\delta$.
\end{Proposition}
\begin{proof}
    Let us prove that $\SpaceV_\SubFamilly$ is an operad ideal of $\SpaceCliff_\delta$. For
    this, let $\BasisF_u, \BasisF_v \in \SpaceCliff_\delta$ and $i \in [|u|]$, and set $f :=
    \BasisF_u \Composition_i \BasisF_v$. We rely on the expression provided by
    Proposition~\ref{prop:composition_f_basis} to compute the partial composition of two
    elements of the $\BasisF$-basis in $\SpaceCliff_\delta$.
    \begin{enumerate}[label=(\arabic*)]
        \item Assume by contradiction that $u \notin \SubFamilly$ and that there is a $w \in
        \SubFamilly$ such that $\BasisF_w$ appears in $f$.  By
        Lemma~\ref{lem:interval_white_and_black_squares}, $\Reduction_\delta\Par{w(1, i - 1)
        \ w(i + \Length(v), \Length(w))} = u$. Since $\SubFamilly$ is closed by subword
        reduction, this would imply that $u \in \SubFamilly$, which contradicts our
        hypothesis.
        \item Similarly, assume now by contradiction that $v \notin \SubFamilly$ and that
        there is a $w \in \SubFamilly$ such that $\BasisF_w$ appears in $f$. By
        Lemma~\ref{lem:interval_white_and_black_squares}, $\Reduction_\delta\Par{w(i,
        \Length(v) - 1)} = v$. Since $\SubFamilly$ is closed by subword reduction, this
        would imply that $v \in \SubFamilly$, which contradicts our hypothesis.
    \end{enumerate}
    Therefore, $f$ belongs in both cases to $\SpaceV_\SubFamilly$. Since moreover $\epsilon
    \in \SubFamilly$, $\BasisF_\epsilon \notin \SpaceV_\SubFamilly$, implying that
    $\SpaceV_\SubFamilly$ is an operad ideal of $\SpaceCliff_\delta$.
\end{proof}

\subsection{Partial composition maps}
We provide now expressions to compute the partial composition maps on different bases for
the quotient $\SpaceCliff_\SubFamilly$ of $\SpaceCliff_\delta$.

We shall us in the sequel the canonical projection map $\theta_\SubFamilly :
\SpaceCliff_\delta \to \SpaceCliff_\SubFamilly$ satisfying, for any $w \in
\SetCliff_\delta$,
\begin{equation}
    \theta_\SubFamilly\Par{\BasisF_w} =
    \begin{cases}
        \BasisF_w & \mbox{if } w \in \SubFamilly, \\
        0 & \mbox{otherwise}.
    \end{cases}
\end{equation}

\subsubsection{Over the fundamental basis}
We first need a property of $\SubFamilly$ depending upon the fact that each
$\SubFamilly(n)$, $n \geq 1$, forms a sublattice of the $\delta$-cliff poset.

\begin{Lemma} \label{lem:subfamilly_sublattice_meet_join}
    Let $\delta$ be a range map and $\SubFamilly$ be a graded subset of $\SetCliff_\delta$
    such that for any $n \geq 1$, $\SubFamilly(n)$ is a sublattice of $\SetCliff_\delta(n)$.
    For any $w \in \SetCliff_\delta(n)$,
    \begin{enumerate}[label={\it (\roman*)}]
        \item \label{item:subfamilly_sublattice_meet_join_1}
        the set $\Bra{w' \in \SubFamilly : w \Leq w'}$ admits at most one minimal element;
        \item \label{item:subfamilly_sublattice_meet_join_2}
        the set $\Bra{w' \in \SubFamilly : w' \Leq w}$ admits at most one maximal element.
    \end{enumerate}
\end{Lemma}
\begin{proof}
    Let $u, v \in X$, where $X$ is the set considered
    in~\ref{item:subfamilly_sublattice_meet_join_1}. Since $\SubFamilly(n)$ is a sublattice
    of $\SetCliff_\delta(n)$, the meet $w'$ of $u$ and $v$ is an element of $\SubFamilly$.
    Moreover, by definition of $X$, $w$ is a lower bound of $\{u, v\}$. Therefore, since
    $w'$ is the greatest lower bound of $\{u, v\}$, we have $w \Leq w'$. For this reason,
    $w' \in X$. This implies~\ref{item:subfamilly_sublattice_meet_join_1}. Similar arguments
    show~\ref{item:subfamilly_sublattice_meet_join_2}.
\end{proof}

As a consequence of Lemma~\ref{lem:subfamilly_sublattice_meet_join}, when $\SubFamilly$
satisfies the given prerequisites, for any $u, v \in \SetCliff_\delta(n)$ such that $u \Leq
v$, $\Han{u, v}_\Leq \cap \SubFamilly$ is empty or is an interval of $\SubFamilly(n)$.
Moreover, let us denote by $\Meet_\SubFamilly(w)$ (resp.\ $\JJoin_\SubFamilly(w)$) the
unique minimal (resp.\ maximal) element of the set described
in~\ref{item:subfamilly_sublattice_meet_join_1} (resp.\
\ref{item:subfamilly_sublattice_meet_join_2}) when it is nonempty.

\begin{Theorem} \label{thm:composition_f_basis_quotient}
    Let $\delta$ be a unimodal range map and $\SubFamilly$ be a nonempty graded subset of
    $\SetCliff_\delta$ such that $\SubFamilly$ is closed by subword reduction. For any $u, v
    \in \SubFamilly$ and $i \in [|u|]$,
    \begin{equation} \label{equ:composition_f_basis_quotient}
        \BasisF_u \Composition_i \BasisF_v
        = \CharacteristicCliff_\delta\Par{u \WhiteSquare_i v}
        \sum_{w \in \Han{u \WhiteSquare_i v,u \BlackSquare_i v}_\Leq \cap \SubFamilly}
        \BasisF_w.
    \end{equation}
    Moreover, when for any $n \geq 1$, $\SubFamilly(n)$ is a sublattice of
    $\SetCliff_\delta(n)$, if~\eqref{equ:composition_f_basis_quotient} is different from
    $0$, the support of this element is the interval $\Han{\Meet_\SubFamilly\Par{u
    \WhiteSquare_i v}, \JJoin_\SubFamilly\Par{u \BlackSquare_i v}}_\Leq$ of $\SubFamilly$.
\end{Theorem}
\begin{proof}
    Since by Proposition~\ref{prop:quotient_cliff_operads}, $\SpaceCliff_\SubFamilly$ is a
    quotient operad of $\SpaceCliff_\delta$,
    \begin{math}
        \BasisF_u \Composition_i \BasisF_v
        = \theta_\SubFamilly\Par{\BasisF_u \Composition_i \BasisF_v},
    \end{math}
    where the partial composition in the left-hand side (resp.\ right-hand side) is the one
    of $\SpaceCliff_\SubFamilly$ (resp.\ $\SpaceCliff_\delta$).
    Expression~\eqref{equ:composition_f_basis_quotient} follows now from
    Proposition~\ref{prop:composition_f_basis} and the linearity of the canonical projection
    map~$\theta_\SubFamilly$.

    The second part of the statement of theorem is implied by the first part and by
    Lemma~\ref{lem:subfamilly_sublattice_meet_join}.
\end{proof}

\subsubsection{Over the elementary basis}
For any $w \in \SubFamilly$, let
\begin{equation} \label{equ:e_basis_from_f_basis_quotient}
    \BasisE_w := \theta_\SubFamilly\Par{\BasisE_w}
    = \sum_{\substack{
        w' \in \SubFamilly \\
        w \Leq w'
    }}
    \BasisF_{w'},
\end{equation}
where the second occurrence of $\BasisE_w$ is an element of $\SpaceCliff_\delta$. By
triangularity, the set $\Bra{\BasisE_w : w \in \SubFamilly}$ is a basis of
$\SpaceCliff_\SubFamilly$, called \Def{elementary basis} (or \Def{$\BasisE$-basis} for
short).

\begin{Proposition} \label{prop:composition_e_basis_quotient}
    Let $\delta$ be a unimodal range map and $\SubFamilly$ be a nonempty graded subset of
    $\SetCliff_\delta$ such that $\SubFamilly$ is closed by subword reduction, and for any
    $n \geq 1$, $\SubFamilly(n)$ is a sublattice of $\SetCliff_\delta(n)$. For any $u, v
    \in \SubFamilly$ and $i \in [|u|]$,
    \begin{equation} \label{equ:composition_e_basis_quotient}
        \BasisE_u \Composition_i \BasisE_v =
        \begin{cases}
            \CharacteristicCliff_\delta\Par{u \WhiteSquare_i v}
            \BasisE_{\Meet_\SubFamilly\Par{u \WhiteSquare_i v}}
                & \mbox{if }
                \Bra{w \in \SubFamilly : u \WhiteSquare_i v \Leq w} \ne \emptyset, \\
            0 & \mbox{otherwise}.
        \end{cases}
    \end{equation}
\end{Proposition}
\begin{proof}
    Since $\theta_\SubFamilly$ is an operad morphism,
    by~\eqref{equ:e_basis_from_f_basis_quotient} and~\eqref{equ:composition_e_basis},
    \begin{equation} \begin{split} \label{equ:composition_e_basis_quotient_1}
        \BasisE_u \Composition_i \BasisE_v
        & =
        \theta_\SubFamilly\Par{\BasisE_u} \Composition_i \theta_\SubFamilly\Par{\BasisE_v}
        =
        \theta_\SubFamilly\Par{\BasisE_u \Composition_i \BasisE_v} \\
        & =
        \theta_\SubFamilly\Par{
        \CharacteristicCliff_\delta\Par{u \WhiteSquare_i v} \BasisE_{u \WhiteSquare_i v}} \\
        & =
        \CharacteristicCliff_\delta\Par{u \WhiteSquare_i v}
        \sum_{\substack{
            w \in \SubFamilly \\
            u \WhiteSquare_i v \Leq w
        }}
        \BasisF_w.
    \end{split} \end{equation}
    By Lemma~\ref{lem:subfamilly_sublattice_meet_join}, when $u \WhiteSquare_i v$ is a
    $\delta$-cliff and $\Bra{w \in \SubFamilly : u \WhiteSquare_i v \Leq w} \ne \emptyset$,
    $\Meet_\SubFamilly\Par{u \WhiteSquare_i v}$ is a well-defined element of $\SubFamilly$.
    In this case, the last term of~\eqref{equ:composition_e_basis_quotient_1} is equal to
    $\BasisE_{\Meet_\SubFamilly\Par{u \WhiteSquare_i v}}$. Otherwise, when $\Bra{w \in
    \SubFamilly : u \WhiteSquare_i v \Leq w} = \emptyset$, the last term
    of~\eqref{equ:composition_e_basis_quotient_1} is zero. This
    establishes~\eqref{equ:composition_e_basis_quotient}.
\end{proof}

In the forthcoming Section~\ref{sec:particular_cases}, we shall consider minimal generating
sets and nontrivial relations of $\SpaceCliff_\SubFamilly$ expressed on the $\BasisE$-basis.
For this reason, we denote by $\GeneratingSet_\SubFamilly$ the unique minimal generating set
of $\SpaceCliff_\SubFamilly$ which is a subset of the $\BasisE$-basis, and by
$\RelationSpace_\SubFamilly$ the space of the nontrivial relations of
$\SpaceCliff_\SubFamilly$.

\subsubsection{Over the homogeneous basis}
For any $w \in \SubFamilly$, let
\begin{equation} \label{equ:h_basis_from_f_basis_quotient}
    \BasisH_w := \theta_\SubFamilly\Par{\BasisH_w}
    = \sum_{\substack{
        w' \in \SubFamilly \\
        w' \Leq w
    }}
    \BasisF_{w'},
\end{equation}
where the second occurrence of $\BasisH_w$ is an element of $\SpaceCliff_\delta$. By
triangularity, the set $\Bra{\BasisH_w : w \in \SubFamilly}$ is a basis of
$\SpaceCliff_\SubFamilly$, called \Def{homogeneous basis} (or \Def{$\BasisH$-basis} for
short).

\begin{Proposition} \label{prop:composition_h_basis_quotient}
    Let $\delta$ be a unimodal range map and $\SubFamilly$ be a nonempty graded subset of
    $\SetCliff_\delta$ such that $\SubFamilly$ is closed by subword reduction, and for any
    $n \geq 1$, $\SubFamilly(n)$ is a sublattice of $\SetCliff_\delta(n)$. For any $u, v
    \in \SubFamilly$ and $i \in [|u|]$,
    \begin{equation} \label{equ:composition_h_basis_quotient}
        \BasisH_u \Composition_i \BasisH_v =
        \begin{cases}
            \BasisH_{\JJoin_\SubFamilly\Par{\Reduction_\delta\Par{u \BlackSquare_i v}}}
                & \mbox{if }
                \Bra{w \in \SubFamilly : w \Leq u \BlackSquare_i v} \ne \emptyset, \\
            0 & \mbox{otherwise}.
        \end{cases}
    \end{equation}
\end{Proposition}
\begin{proof}
    Since $\theta_\SubFamilly$ is an operad morphism,
    by~\eqref{equ:h_basis_from_f_basis_quotient} and
    Proposition~\ref{prop:composition_h_basis},
    \begin{equation} \begin{split} \label{equ:composition_h_basis_quotient_1}
        \BasisH_u \Composition_i \BasisH_v
        & =
        \theta_\SubFamilly\Par{\BasisH_u} \Composition_i \theta_\SubFamilly\Par{\BasisH_v}
        =
        \theta_\SubFamilly\Par{\BasisH_u \Composition_i \BasisH_v} \\
        & =
        \theta_\SubFamilly\Par{\BasisH_{\Reduction_\delta\Par{u \BlackSquare_i v}}} \\
        & =
        \sum_{\substack{
            w \in \SubFamilly \\
            w \Leq \Reduction_\delta\Par{u \BlackSquare_i v}
        }}
        \BasisF_w.
    \end{split} \end{equation}
    By Lemma~\ref{lem:subfamilly_sublattice_meet_join}, when $\Bra{w \in \SubFamilly : w
    \Leq u \BlackSquare_i v} \ne \emptyset$, $\JJoin_\SubFamilly\Par{u \BlackSquare_i v}$ is
    a well-defined element of $\SubFamilly$. In this case, the last term
    of~\eqref{equ:composition_h_basis_quotient_1} is equal to
    $\BasisH_{\JJoin_\SubFamilly\Par{\Reduction_\delta\Par{u \BlackSquare_i v}}}$.
    Otherwise, when $\Bra{w \in \SubFamilly : w \Leq u \BlackSquare_i v} = \emptyset$, the
    last term of~\eqref{equ:composition_h_basis_quotient_1} is zero. This
    establishes~\eqref{equ:composition_h_basis_quotient}.
\end{proof}

\section{Some particular constructions} \label{sec:particular_cases}
We study in this last part some operads $\SpaceCliff_\delta$ and quotients
$\SpaceCliff_\SubFamilly$ for some concrete range maps and graded sets $\SubFamilly$ of
$\delta$-cliffs.

\subsection{Operads on cliffs}
Let us begin by providing some properties about the operads $\SpaceCliff_\MapConstant{c}$,
$c \geq 0$, and $\SpaceCliff_\MapArithmetic{m}$, $m \geq 0$.

\subsubsection{On constant range maps}
For any $c \geq 0$, the operad $\SpaceCliff_\MapConstant{c}$ is by construction the
interstice operad $\IntersticeOperad(\HanL{c})$. Therefore, $\SpaceCliff_\MapConstant{c}$
has the properties presented in Section~\ref{subsubsec:interstice_operads}. As a particular
case, the map $\ToComposition$ (see Section~\ref{subsubsec:integer_compositions}) allows us
to interpret any $\MapConstant{1}$-cliff as an integer composition. Therefore,
$\SpaceCliff_{\MapConstant{1}}$ can be seen as an operad on integers compositions. For
instance,
\begin{subequations}
\begin{equation}
    \BasisE_{(1, 2, 1, 2, 2)} \Composition_5 \BasisE_{(2, 3, 1, 1)}
    = \BasisE_{(1, 2, 1, 2, 3, 1, 2, 2)},
\end{equation}
\begin{equation}
    \BasisF_{(1, 2, 1, 2, 2)} \Composition_5 \BasisF_{(2, 3, 1, 1)}
    = \BasisF_{(1, 2, 1, 2, 3, 1, 2, 2)},
\end{equation}
\begin{equation}
    \BasisH_{(1, 2, 1, 2, 2)} \Composition_5 \BasisH_{(2, 3, 1, 1)}
    = \BasisH_{(1, 2, 1, 2, 3, 1, 2, 2)}.
\end{equation}
\end{subequations}
It is possible to prove, by using Propositions~\ref{prop:composition_f_basis}
and~\ref{prop:composition_h_basis}, that the constant structures of
$\SpaceCliff_{\MapConstant{1}}$ are the same in the $\BasisE$, $\BasisF$, and
$\BasisH$-bases.

\subsubsection{On arithmetic range maps}
Let us study the operads $\SpaceCliff_\MapArithmetic{m}$, $m \geq 0$. The map
$\ToPermutation$ (see Section~\ref{subsubsec:permutations}) allows us to interpret any
$\MapArithmetic{1}$-cliff as a permutation. Therefore, $\SpaceCliff_{\MapArithmetic{1}}$ can
be seen as an operad on permutations. For instance,
\begin{subequations}
\begin{equation}
    \BasisE_{25143} \Composition_3 \BasisE_{3142} = \BasisE_{215369487},
\end{equation}
\begin{multline}
    \BasisF_{25143} \Composition_3 \BasisF_{3142}
    = \BasisF_{215369487} + \BasisF_{251369487} + \BasisF_{521369487} + \BasisF_{235169487}
    + \BasisF_{253169487} + \BasisF_{523169487} \\
    + \BasisF_{325196487} + \BasisF_{352169487} + \BasisF_{532169487},
\end{multline}
\begin{equation}
    \BasisH_{25143} \Composition_3 \BasisH_{3142} = \BasisH_{532169487}.
\end{equation}
\end{subequations}
In the same way, the map $\ToTree$ (see Section~\ref{subsubsec:increasing_trees}) allows us
to interpret any $\MapArithmetic{m}$-cliff as an $m$-increasing tree. Therefore,
$\SpaceCliff_{\MapArithmetic{m}}$ can be seen as an operad on $m$-increasing trees. For
instance, in $\SpaceCliff_\MapArithmetic{2}$,
\begin{subequations}
\begin{equation}
    \BasisE_{
        \scalebox{.65}{
        \begin{tikzpicture}[Centering,xscale=0.17,yscale=0.14]
            \node[Leaf](0)at(1.00,-3.25){};
            \node[Leaf](11)at(7.00,-6.50){};
            \node[Leaf](12)at(8.00,-6.50){};
            \node[Leaf](2)at(1.00,-9.75){};
            \node[Leaf](4)at(2.00,-9.75){};
            \node[Leaf](5)at(3.00,-9.75){};
            \node[Leaf](7)at(4.00,-6.50){};
            \node[Leaf](8)at(5.00,-6.50){};
            \node[Leaf](9)at(6.00,-6.50){};
            \node[Node](1)at(4.00,0.00){$1$};
            \node[Node](10)at(7.00,-3.25){$2$};
            \node[Node](3)at(2.00,-6.50){$4$};
            \node[Node](6)at(4.00,-3.25){$3$};
            \draw[Edge](0)--(1);
            \draw[Edge](10)--(1);
            \draw[Edge](11)--(10);
            \draw[Edge](12)--(10);
            \draw[Edge](2)--(3);
            \draw[Edge](3)--(6);
            \draw[Edge](4)--(3);
            \draw[Edge](5)--(3);
            \draw[Edge](6)--(1);
            \draw[Edge](7)--(6);
            \draw[Edge](8)--(6);
            \draw[Edge](9)--(10);
            \node(r)at(4.00,3){};
            \draw[Edge](r)--(1);
        \end{tikzpicture}}}
    \Composition_2
    \BasisE_{
        \scalebox{.65}{
        \begin{tikzpicture}[Centering,xscale=0.19,yscale=0.18]
            \node[Leaf](0)at(0.00,-2.50){};
            \node[Leaf](2)at(1.00,-5.00){};
            \node[Leaf](4)at(2.00,-5.00){};
            \node[Leaf](5)at(3.00,-7.50){};
            \node[Leaf](7)at(4.00,-7.50){};
            \node[Leaf](8)at(5.00,-7.50){};
            \node[Leaf](9)at(4.00,-2.50){};
            \node[Node,MarkA](1)at(2.00,0.00){$1$};
            \node[Node,MarkA](3)at(2.00,-2.50){$2$};
            \node[Node,MarkA](6)at(4.00,-5.00){$3$};
            \draw[Edge](0)--(1);
            \draw[Edge](2)--(3);
            \draw[Edge](3)--(1);
            \draw[Edge](4)--(3);
            \draw[Edge](5)--(6);
            \draw[Edge](6)--(3);
            \draw[Edge](7)--(6);
            \draw[Edge](8)--(6);
            \draw[Edge](9)--(1);
            \node(r)at(2.00,2){};
            \draw[Edge](r)--(1);
        \end{tikzpicture}}}
    =
    \BasisE_{
        \scalebox{.65}{
        \begin{tikzpicture}[Centering,xscale=0.17,yscale=0.13]
            \node[Leaf](0)at(3.00,-7.33){};
            \node[Leaf](10)at(6.00,-11.00){};
            \node[Leaf](11)at(7.00,-14.67){};
            \node[Leaf](13)at(8.00,-14.67){};
            \node[Leaf](14)at(9.00,-14.67){};
            \node[Leaf](15)at(9.00,-7.33){};
            \node[Leaf](17)at(9.5,-3.67){};
            \node[Leaf](18)at(12.00,-7.33){};
            \node[Leaf](2)at(1.00,-18.33){};
            \node[Leaf](20)at(13.00,-7.33){};
            \node[Leaf](21)at(14.00,-7.33){};
            \node[Leaf](4)at(2.00,-18.33){};
            \node[Leaf](5)at(3.00,-18.33){};
            \node[Leaf](7)at(4.00,-14.67){};
            \node[Leaf](8)at(5.00,-14.67){};
            \node[Node,MarkA](1)at(6.00,-3.67){$3$};
            \node[Node,MarkA](12)at(8.00,-11.00){$5$};
            \node[Node](16)at(9.5,0.00){$1$};
            \node[Node](19)at(13.00,-3.67){$2$};
            \node[Node](3)at(2.00,-14.67){$7$};
            \node[Node](6)at(4.00,-11.00){$6$};
            \node[Node,MarkA](9)at(6.00,-7.33){$4$};
            \draw[Edge](0)--(1);
            \draw[Edge](1)--(16);
            \draw[Edge](10)--(9);
            \draw[Edge](11)--(12);
            \draw[Edge](12)--(9);
            \draw[Edge](13)--(12);
            \draw[Edge](14)--(12);
            \draw[Edge](15)--(1);
            \draw[Edge](17)--(16);
            \draw[Edge](18)--(19);
            \draw[Edge](19)--(16);
            \draw[Edge](2)--(3);
            \draw[Edge](20)--(19);
            \draw[Edge](21)--(19);
            \draw[Edge](3)--(6);
            \draw[Edge](4)--(3);
            \draw[Edge](5)--(3);
            \draw[Edge](6)--(9);
            \draw[Edge](7)--(6);
            \draw[Edge](8)--(6);
            \draw[Edge](9)--(1);
            \node(r)at(9.5,3){};
            \draw[Edge](r)--(16);
        \end{tikzpicture}}},
\end{equation}
\begin{equation}
    \BasisF_{
        \scalebox{.65}{
        \begin{tikzpicture}[Centering,xscale=0.17,yscale=0.14]
            \node[Leaf](0)at(1.00,-3.25){};
            \node[Leaf](11)at(7.00,-6.50){};
            \node[Leaf](12)at(8.00,-6.50){};
            \node[Leaf](2)at(1.00,-9.75){};
            \node[Leaf](4)at(2.00,-9.75){};
            \node[Leaf](5)at(3.00,-9.75){};
            \node[Leaf](7)at(4.00,-6.50){};
            \node[Leaf](8)at(5.00,-6.50){};
            \node[Leaf](9)at(6.00,-6.50){};
            \node[Node](1)at(4.00,0.00){$1$};
            \node[Node](10)at(7.00,-3.25){$2$};
            \node[Node](3)at(2.00,-6.50){$4$};
            \node[Node](6)at(4.00,-3.25){$3$};
            \draw[Edge](0)--(1);
            \draw[Edge](10)--(1);
            \draw[Edge](11)--(10);
            \draw[Edge](12)--(10);
            \draw[Edge](2)--(3);
            \draw[Edge](3)--(6);
            \draw[Edge](4)--(3);
            \draw[Edge](5)--(3);
            \draw[Edge](6)--(1);
            \draw[Edge](7)--(6);
            \draw[Edge](8)--(6);
            \draw[Edge](9)--(10);
            \node(r)at(4.00,3){};
            \draw[Edge](r)--(1);
        \end{tikzpicture}}}
    \Composition_2
    \BasisF_{
        \scalebox{.65}{
        \begin{tikzpicture}[Centering,xscale=0.19,yscale=0.18]
            \node[Leaf](0)at(0.00,-2.50){};
            \node[Leaf](2)at(1.00,-5.00){};
            \node[Leaf](4)at(2.00,-5.00){};
            \node[Leaf](5)at(3.00,-7.50){};
            \node[Leaf](7)at(4.00,-7.50){};
            \node[Leaf](8)at(5.00,-7.50){};
            \node[Leaf](9)at(4.00,-2.50){};
            \node[Node,MarkA](1)at(2.00,0.00){$1$};
            \node[Node,MarkA](3)at(2.00,-2.50){$2$};
            \node[Node,MarkA](6)at(4.00,-5.00){$3$};
            \draw[Edge](0)--(1);
            \draw[Edge](2)--(3);
            \draw[Edge](3)--(1);
            \draw[Edge](4)--(3);
            \draw[Edge](5)--(6);
            \draw[Edge](6)--(3);
            \draw[Edge](7)--(6);
            \draw[Edge](8)--(6);
            \draw[Edge](9)--(1);
            \node(r)at(2.00,2){};
            \draw[Edge](r)--(1);
        \end{tikzpicture}}}
    =
    \BasisF_{
        \scalebox{.65}{
        \begin{tikzpicture}[Centering,xscale=0.17,yscale=0.13]
            \node[Leaf](0)at(3.00,-7.33){};
            \node[Leaf](10)at(6.00,-11.00){};
            \node[Leaf](11)at(7.00,-14.67){};
            \node[Leaf](13)at(8.00,-14.67){};
            \node[Leaf](14)at(9.00,-14.67){};
            \node[Leaf](15)at(9.00,-7.33){};
            \node[Leaf](17)at(9.5,-3.67){};
            \node[Leaf](18)at(12.00,-7.33){};
            \node[Leaf](2)at(1.00,-18.33){};
            \node[Leaf](20)at(13.00,-7.33){};
            \node[Leaf](21)at(14.00,-7.33){};
            \node[Leaf](4)at(2.00,-18.33){};
            \node[Leaf](5)at(3.00,-18.33){};
            \node[Leaf](7)at(4.00,-14.67){};
            \node[Leaf](8)at(5.00,-14.67){};
            \node[Node,MarkA](1)at(6.00,-3.67){$3$};
            \node[Node,MarkA](12)at(8.00,-11.00){$5$};
            \node[Node](16)at(9.5,0.00){$1$};
            \node[Node](19)at(13.00,-3.67){$2$};
            \node[Node](3)at(2.00,-14.67){$7$};
            \node[Node](6)at(4.00,-11.00){$6$};
            \node[Node,MarkA](9)at(6.00,-7.33){$4$};
            \draw[Edge](0)--(1);
            \draw[Edge](1)--(16);
            \draw[Edge](10)--(9);
            \draw[Edge](11)--(12);
            \draw[Edge](12)--(9);
            \draw[Edge](13)--(12);
            \draw[Edge](14)--(12);
            \draw[Edge](15)--(1);
            \draw[Edge](17)--(16);
            \draw[Edge](18)--(19);
            \draw[Edge](19)--(16);
            \draw[Edge](2)--(3);
            \draw[Edge](20)--(19);
            \draw[Edge](21)--(19);
            \draw[Edge](3)--(6);
            \draw[Edge](4)--(3);
            \draw[Edge](5)--(3);
            \draw[Edge](6)--(9);
            \draw[Edge](7)--(6);
            \draw[Edge](8)--(6);
            \draw[Edge](9)--(1);
            \node(r)at(9.5,3){};
            \draw[Edge](r)--(16);
        \end{tikzpicture}}}
    +
    \BasisF_{
        \scalebox{.65}{
        \begin{tikzpicture}[Centering,xscale=0.17,yscale=0.13]
            \node[Leaf](0)at(7.00,-3.67){};
            \node[Leaf](10)at(6.00,-11.00){};
            \node[Leaf](11)at(7.00,-14.67){};
            \node[Leaf](13)at(8.00,-14.67){};
            \node[Leaf](14)at(9.00,-14.67){};
            \node[Leaf](16)at(10.00,-7.33){};
            \node[Leaf](17)at(11.00,-7.33){};
            \node[Leaf](18)at(12.00,-7.33){};
            \node[Leaf](2)at(1.00,-18.33){};
            \node[Leaf](20)at(13.00,-7.33){};
            \node[Leaf](21)at(14.00,-7.33){};
            \node[Leaf](4)at(2.00,-18.33){};
            \node[Leaf](5)at(3.00,-18.33){};
            \node[Leaf](7)at(4.00,-14.67){};
            \node[Leaf](8)at(5.00,-14.67){};
            \node[Node](1)at(10.00,0.00){$1$};
            \node[Node,MarkA](12)at(8.00,-11.00){$5$};
            \node[Node,MarkA](15)at(10.00,-3.67){$3$};
            \node[Node](19)at(13.00,-3.67){$2$};
            \node[Node](3)at(2.00,-14.67){$7$};
            \node[Node](6)at(4.00,-11.00){$6$};
            \node[Node,MarkA](9)at(6.00,-7.33){$4$};
            \draw[Edge](0)--(1);
            \draw[Edge](10)--(9);
            \draw[Edge](11)--(12);
            \draw[Edge](12)--(9);
            \draw[Edge](13)--(12);
            \draw[Edge](14)--(12);
            \draw[Edge](15)--(1);
            \draw[Edge](16)--(15);
            \draw[Edge](17)--(15);
            \draw[Edge](18)--(19);
            \draw[Edge](19)--(1);
            \draw[Edge](2)--(3);
            \draw[Edge](20)--(19);
            \draw[Edge](21)--(19);
            \draw[Edge](3)--(6);
            \draw[Edge](4)--(3);
            \draw[Edge](5)--(3);
            \draw[Edge](6)--(9);
            \draw[Edge](7)--(6);
            \draw[Edge](8)--(6);
            \draw[Edge](9)--(15);
            \node(r)at(10.00,2.75){};
            \draw[Edge](r)--(1);
        \end{tikzpicture}}}
    +
    \BasisF_{
        \scalebox{.65}{
        \begin{tikzpicture}[Centering,xscale=0.17,yscale=0.13]
            \node[Leaf](0)at(-1.00,-4.40){};
            \node[Leaf](10)at(6.00,-8.80){};
            \node[Leaf](11)at(7.00,-13.20){};
            \node[Leaf](13)at(8.00,-13.20){};
            \node[Leaf](14)at(9.00,-13.20){};
            \node[Leaf](15)at(10.00,-13.20){};
            \node[Leaf](17)at(11.00,-13.20){};
            \node[Leaf](18)at(12.00,-13.20){};
            \node[Leaf](2)at(1.00,-17.60){};
            \node[Leaf](20)at(13.00,-8.80){};
            \node[Leaf](21)at(14.00,-8.80){};
            \node[Leaf](4)at(2.00,-17.60){};
            \node[Leaf](5)at(3.00,-17.60){};
            \node[Leaf](7)at(4.00,-13.20){};
            \node[Leaf](8)at(5.00,-13.20){};
            \node[Node](1)at(6.00,0.00){$1$};
            \node[Node,MarkA](12)at(8.00,-8.80){$5$};
            \node[Node,MarkA](16)at(11.00,-8.80){$3$};
            \node[Node](19)at(13.00,-4.40){$2$};
            \node[Node](3)at(2.00,-13.20){$7$};
            \node[Node](6)at(4.00,-8.80){$6$};
            \node[Node,MarkA](9)at(6.00,-4.40){$4$};
            \draw[Edge](0)--(1);
            \draw[Edge](10)--(9);
            \draw[Edge](11)--(12);
            \draw[Edge](12)--(9);
            \draw[Edge](13)--(12);
            \draw[Edge](14)--(12);
            \draw[Edge](15)--(16);
            \draw[Edge](16)--(19);
            \draw[Edge](17)--(16);
            \draw[Edge](18)--(16);
            \draw[Edge](19)--(1);
            \draw[Edge](2)--(3);
            \draw[Edge](20)--(19);
            \draw[Edge](21)--(19);
            \draw[Edge](3)--(6);
            \draw[Edge](4)--(3);
            \draw[Edge](5)--(3);
            \draw[Edge](6)--(9);
            \draw[Edge](7)--(6);
            \draw[Edge](8)--(6);
            \draw[Edge](9)--(1);
            \node(r)at(6.00,3.30){};
            \draw[Edge](r)--(1);
        \end{tikzpicture}}}
    +
    \BasisF_{
        \scalebox{.65}{
        \begin{tikzpicture}[Centering,xscale=0.17,yscale=0.13]
            \node[Leaf](0)at(0.00,-4.40){};
            \node[Leaf](10)at(6.00,-8.80){};
            \node[Leaf](11)at(7.00,-13.20){};
            \node[Leaf](13)at(8.00,-13.20){};
            \node[Leaf](14)at(9.00,-13.20){};
            \node[Leaf](15)at(10.00,-8.80){};
            \node[Leaf](17)at(11.00,-13.20){};
            \node[Leaf](19)at(12.00,-13.20){};
            \node[Leaf](2)at(1.00,-17.60){};
            \node[Leaf](20)at(13.00,-13.20){};
            \node[Leaf](21)at(14.00,-8.80){};
            \node[Leaf](4)at(2.00,-17.60){};
            \node[Leaf](5)at(3.00,-17.60){};
            \node[Leaf](7)at(4.00,-13.20){};
            \node[Leaf](8)at(5.00,-13.20){};
            \node[Node](1)at(6.00,0.00){$1$};
            \node[Node,MarkA](12)at(8.00,-8.80){$5$};
            \node[Node](16)at(12.00,-4.40){$2$};
            \node[Node,MarkA](18)at(12.00,-8.80){$3$};
            \node[Node](3)at(2.00,-13.20){$7$};
            \node[Node](6)at(4.00,-8.80){$6$};
            \node[Node,MarkA](9)at(6.00,-4.40){$4$};
            \draw[Edge](0)--(1);
            \draw[Edge](10)--(9);
            \draw[Edge](11)--(12);
            \draw[Edge](12)--(9);
            \draw[Edge](13)--(12);
            \draw[Edge](14)--(12);
            \draw[Edge](15)--(16);
            \draw[Edge](16)--(1);
            \draw[Edge](17)--(18);
            \draw[Edge](18)--(16);
            \draw[Edge](19)--(18);
            \draw[Edge](2)--(3);
            \draw[Edge](20)--(18);
            \draw[Edge](21)--(16);
            \draw[Edge](3)--(6);
            \draw[Edge](4)--(3);
            \draw[Edge](5)--(3);
            \draw[Edge](6)--(9);
            \draw[Edge](7)--(6);
            \draw[Edge](8)--(6);
            \draw[Edge](9)--(1);
            \node(r)at(6.00,3.30){};
            \draw[Edge](r)--(1);
        \end{tikzpicture}}}
    +
    \BasisF_{
        \scalebox{.65}{
        \begin{tikzpicture}[Centering,xscale=0.17,yscale=0.12]
            \node[Leaf](0)at(1.00,-4.40){};
            \node[Leaf](10)at(6.00,-8.80){};
            \node[Leaf](11)at(7.00,-13.20){};
            \node[Leaf](13)at(8.00,-13.20){};
            \node[Leaf](14)at(9.00,-13.20){};
            \node[Leaf](15)at(10.00,-8.80){};
            \node[Leaf](17)at(11.00,-8.80){};
            \node[Leaf](18)at(12.00,-13.20){};
            \node[Leaf](2)at(1.00,-17.60){};
            \node[Leaf](20)at(13.00,-13.20){};
            \node[Leaf](21)at(14.00,-13.20){};
            \node[Leaf](4)at(2.00,-17.60){};
            \node[Leaf](5)at(3.00,-17.60){};
            \node[Leaf](7)at(4.00,-13.20){};
            \node[Leaf](8)at(5.00,-13.20){};
            \node[Node](1)at(6.00,0.00){$1$};
            \node[Node,MarkA](12)at(8.00,-8.80){$5$};
            \node[Node](16)at(11.00,-4.40){$2$};
            \node[Node,MarkA](19)at(13.00,-8.80){$3$};
            \node[Node](3)at(2.00,-13.20){$7$};
            \node[Node](6)at(4.00,-8.80){$6$};
            \node[Node,MarkA](9)at(6.00,-4.40){$4$};
            \draw[Edge](0)--(1);
            \draw[Edge](10)--(9);
            \draw[Edge](11)--(12);
            \draw[Edge](12)--(9);
            \draw[Edge](13)--(12);
            \draw[Edge](14)--(12);
            \draw[Edge](15)--(16);
            \draw[Edge](16)--(1);
            \draw[Edge](17)--(16);
            \draw[Edge](18)--(19);
            \draw[Edge](19)--(16);
            \draw[Edge](2)--(3);
            \draw[Edge](20)--(19);
            \draw[Edge](21)--(19);
            \draw[Edge](3)--(6);
            \draw[Edge](4)--(3);
            \draw[Edge](5)--(3);
            \draw[Edge](6)--(9);
            \draw[Edge](7)--(6);
            \draw[Edge](8)--(6);
            \draw[Edge](9)--(1);
            \node(r)at(6.00,3.30){};
            \draw[Edge](r)--(1);
        \end{tikzpicture}}},
\end{equation}
\begin{equation}
    \BasisH_{
        \scalebox{.65}{
        \begin{tikzpicture}[Centering,xscale=0.17,yscale=0.14]
            \node[Leaf](0)at(1.00,-3.25){};
            \node[Leaf](11)at(7.00,-6.50){};
            \node[Leaf](12)at(8.00,-6.50){};
            \node[Leaf](2)at(1.00,-9.75){};
            \node[Leaf](4)at(2.00,-9.75){};
            \node[Leaf](5)at(3.00,-9.75){};
            \node[Leaf](7)at(4.00,-6.50){};
            \node[Leaf](8)at(5.00,-6.50){};
            \node[Leaf](9)at(6.00,-6.50){};
            \node[Node](1)at(4.00,0.00){$1$};
            \node[Node](10)at(7.00,-3.25){$2$};
            \node[Node](3)at(2.00,-6.50){$4$};
            \node[Node](6)at(4.00,-3.25){$3$};
            \draw[Edge](0)--(1);
            \draw[Edge](10)--(1);
            \draw[Edge](11)--(10);
            \draw[Edge](12)--(10);
            \draw[Edge](2)--(3);
            \draw[Edge](3)--(6);
            \draw[Edge](4)--(3);
            \draw[Edge](5)--(3);
            \draw[Edge](6)--(1);
            \draw[Edge](7)--(6);
            \draw[Edge](8)--(6);
            \draw[Edge](9)--(10);
            \node(r)at(4.00,3){};
            \draw[Edge](r)--(1);
        \end{tikzpicture}}}
    \Composition_2
    \BasisH_{
        \scalebox{.65}{
        \begin{tikzpicture}[Centering,xscale=0.19,yscale=0.18]
            \node[Leaf](0)at(0.00,-2.50){};
            \node[Leaf](2)at(1.00,-5.00){};
            \node[Leaf](4)at(2.00,-5.00){};
            \node[Leaf](5)at(3.00,-7.50){};
            \node[Leaf](7)at(4.00,-7.50){};
            \node[Leaf](8)at(5.00,-7.50){};
            \node[Leaf](9)at(4.00,-2.50){};
            \node[Node,MarkA](1)at(2.00,0.00){$1$};
            \node[Node,MarkA](3)at(2.00,-2.50){$2$};
            \node[Node,MarkA](6)at(4.00,-5.00){$3$};
            \draw[Edge](0)--(1);
            \draw[Edge](2)--(3);
            \draw[Edge](3)--(1);
            \draw[Edge](4)--(3);
            \draw[Edge](5)--(6);
            \draw[Edge](6)--(3);
            \draw[Edge](7)--(6);
            \draw[Edge](8)--(6);
            \draw[Edge](9)--(1);
            \node(r)at(2.00,2){};
            \draw[Edge](r)--(1);
        \end{tikzpicture}}}
    =
    \BasisH_{
        \scalebox{.65}{
        \begin{tikzpicture}[Centering,xscale=0.17,yscale=0.12]
            \node[Leaf](0)at(1.00,-4.40){};
            \node[Leaf](10)at(6.00,-8.80){};
            \node[Leaf](11)at(7.00,-13.20){};
            \node[Leaf](13)at(8.00,-13.20){};
            \node[Leaf](14)at(9.00,-13.20){};
            \node[Leaf](15)at(10.00,-8.80){};
            \node[Leaf](17)at(11.00,-8.80){};
            \node[Leaf](18)at(12.00,-13.20){};
            \node[Leaf](2)at(1.00,-17.60){};
            \node[Leaf](20)at(13.00,-13.20){};
            \node[Leaf](21)at(14.00,-13.20){};
            \node[Leaf](4)at(2.00,-17.60){};
            \node[Leaf](5)at(3.00,-17.60){};
            \node[Leaf](7)at(4.00,-13.20){};
            \node[Leaf](8)at(5.00,-13.20){};
            \node[Node](1)at(6.00,0.00){$1$};
            \node[Node,MarkA](12)at(8.00,-8.80){$5$};
            \node[Node](16)at(11.00,-4.40){$2$};
            \node[Node,MarkA](19)at(13.00,-8.80){$3$};
            \node[Node](3)at(2.00,-13.20){$7$};
            \node[Node](6)at(4.00,-8.80){$6$};
            \node[Node,MarkA](9)at(6.00,-4.40){$4$};
            \draw[Edge](0)--(1);
            \draw[Edge](10)--(9);
            \draw[Edge](11)--(12);
            \draw[Edge](12)--(9);
            \draw[Edge](13)--(12);
            \draw[Edge](14)--(12);
            \draw[Edge](15)--(16);
            \draw[Edge](16)--(1);
            \draw[Edge](17)--(16);
            \draw[Edge](18)--(19);
            \draw[Edge](19)--(16);
            \draw[Edge](2)--(3);
            \draw[Edge](20)--(19);
            \draw[Edge](21)--(19);
            \draw[Edge](3)--(6);
            \draw[Edge](4)--(3);
            \draw[Edge](5)--(3);
            \draw[Edge](6)--(9);
            \draw[Edge](7)--(6);
            \draw[Edge](8)--(6);
            \draw[Edge](9)--(1);
            \node(r)at(6.00,3.30){};
            \draw[Edge](r)--(1);
        \end{tikzpicture}}}.
\end{equation}
\end{subequations}

Proposition~\ref{prop:minimal_generating_set_cliffs} allows us to list the first elements of
the minimal generating sets $\GeneratingSet_\MapArithmetic{m}$, $m \geq 0$, of
$\SpaceCliff_\MapArithmetic{m}$. Here are the lists of these generators for $m \in
\HanL{2}$, up to arity~$4$:
\begin{subequations}
\begin{equation}
    \BasisE_0,
    \qquad m = 0,
\end{equation}
\begin{equation}
    \BasisE_0, \quad
    \BasisE_{01}, \quad
    \BasisE_{002}, \BasisE_{011}, \BasisE_{012},
    \qquad m = 1,
\end{equation}
\begin{equation}
    \BasisE_0, \quad
    \BasisE_{01}, \BasisE_{02}, \quad
    \BasisE_{003}, \BasisE_{004}, \BasisE_{011}, \BasisE_{012}, \BasisE_{013},
    \BasisE_{014}, \BasisE_{021}, \BasisE_{022}, \BasisE_{023}, \BasisE_{024},
    \qquad m = 2.
\end{equation}
\end{subequations}
By Proposition~\ref{prop:properties_prime_sets}, $\GeneratingSet_\MapArithmetic{0}$ is
finite while $\GeneratingSet_\MapArithmetic{1}$ and $\GeneratingSet_\MapArithmetic{2}$ are
infinite.

In order to enumerate $\GeneratingSet_\MapArithmetic{m}$, we need the following small result
leading to a recursive description of $\MapArithmetic{m}$-prime $\MapArithmetic{m}$-cliffs.

\begin{Lemma} \label{lem:decomposition_zero}
    For any $m \geq 0$, if $w$ is a nonempty $\MapArithmetic{m}$-cliff which is not
    $\MapArithmetic{m}$-prime, then there exists an $\MapArithmetic{m}$-cliff $w'$ and $i
    \in \Han{\Brr{w'}}$ such that $w = w' \WhiteSquare_i 0$.
\end{Lemma}
\begin{proof}
    Since $w$ is not $\MapArithmetic{m}$-prime, there exist two nonempty
    $\MapArithmetic{m}$-cliffs $u$ and $v$, and $i \in [|u|]$ such that $w = u
    \WhiteSquare_i v$. If $i = 1$, then one has the decomposition $w = v \WhiteSquare_{|v|}
    u$ so that, since $|v| \geq 2$, we can assume that $i \geq 2$. Now, since $v$ is an
    $\MapArithmetic{m}$-cliff, we have $v(1) = 0$ and
    \begin{equation} \begin{split}
        w & = u \WhiteSquare_i v \\
        & = u(1, i - 1) \ v \ u(i, \Length(u)) \\
        & = u(1, i - 1) \ 0 \ v(2, \Length(v)) \ u(i, \Length(u)) \\
        & = u(1, i - 1) \ v(2, \Length(v)) \ u(i, \Length(u)) \WhiteSquare_i 0.
    \end{split} \end{equation}
    The fact that $u$ and $v$ are $\MapArithmetic{m}$-cliffs and $i \geq 2$ implies that
    $u(1, i - 1) \ v(2, \Length(v)) \ u(i, \Length(u))$ also is. This establishes the stated
    property.
\end{proof}

\begin{Lemma} \label{lem:prime_m_cliffs}
    Let, for an $m \geq 0$, $w$ be an $\MapArithmetic{m}$-cliff decomposing as $w = w' a$
    with $w' \in \SetCliff_\MapArithmetic{m}$ and $a \in \N$. Then, $w$ is
    $\MapArithmetic{m}$-prime if and only if one of the following two assertions is
    satisfied:
    \begin{enumerate}[label={\it (\roman*)}]
        \item \label{item:prime_m_cliffs_1}
        $w' \notin \SetPrime_\MapArithmetic{m}$ and $a \geq (\Length(w) - 2)m + 1$;
        \item \label{item:prime_m_cliffs_2}
        $w' \in \SetPrime_\MapArithmetic{m}$ and $a \ne 0$.
    \end{enumerate}
\end{Lemma}
\begin{proof}
    Assume first that~\ref{item:prime_m_cliffs_1} holds and assume that there are $u, v \in
    \SetCliff_\MapArithmetic{m}$ and $i \in [|u|]$ such that $w = u \WhiteSquare_i v$. If
    $u$ and $v$ are both nonempty, the letter $a$ of $w$ appears either in $u$ or in $v$ but
    at a position smaller than the one it has in $w$. Since $a \geq (\Length(w) - 2)m + 1$,
    either $u$ of $v$ would not be an $\MapArithmetic{m}$-cliff. Therefore, $u = v =
    \epsilon$ and $w$ is prime. Assume now that~\ref{item:prime_m_cliffs_2} holds and,
    again, assume that there are $u, v \in \SetCliff_\MapArithmetic{m}$ and $i \in [|u|]$
    such that $w = u \WhiteSquare_i v$. If the letter $a$ of $w$ is in $u$, we have $u = u'
    a$ where $u' \in \SetCliff_\MapArithmetic{m}$ and $w' a = u' a \WhiteSquare_i v$.
    Otherwise, the letter $a$ of $w$ is in $v$ and we have $v = v' a$ where $v' \in
    \SetCliff_\MapArithmetic{m}$ and $w' a = u \WhiteSquare_i v'a$. Therefore, $w'$
    decomposes respectively as $w' = u' \WhiteSquare_i v$ and $w' = u \WhiteSquare_i v'$.
    Since $w'$ is $\MapArithmetic{m}$-prime, these decompositions are trivial. This implies
    that $u' = \epsilon$ so that $u = a$, or that $v = \epsilon$, or that $u = \epsilon$, or
    that $v' = \epsilon$ so that $v = a$. Since $a \ne 0$, $a$ is not a
    $\MapArithmetic{m}$-cliff. Therefore, $v = \epsilon$ or $u = \epsilon$, implying that
    $w$ is $\MapArithmetic{m}$-prime.

    Conversely, assume that the negations of~\ref{item:prime_m_cliffs_1}
    and~\ref{item:prime_m_cliffs_2} hold at the same time. Therefore, at least one of the
    following assertions holds.
    \begin{enumerate}[label=(A\arabic*)]
        \item \label{item:prime_m_cliffs_a}
        $w' \in \SetPrime_\MapArithmetic{m}$ and $w' \notin \SetPrime_\MapArithmetic{m}$;
        \item \label{item:prime_m_cliffs_b}
        $w' \in \SetPrime_\MapArithmetic{m}$ and $a = 0$;
        \item \label{item:prime_m_cliffs_c}
        $a < (\Length(w) - 2)m + 1$ and $w' \notin \SetPrime_\MapArithmetic{m}$;
        \item \label{item:prime_m_cliffs_d}
        $a < (\Length(w) - 2)m + 1$ and $a = 0$.
    \end{enumerate}
    Assertion~\ref{item:prime_m_cliffs_a} is absurd so that this situation cannot occurs.
    If~\ref{item:prime_m_cliffs_b} or~\ref{item:prime_m_cliffs_d} holds, then $a = 0$ and
    $w$ decomposes as $w = w' \WhiteSquare_{\Brr{w'}} 0$, showing that $w$ is not
    $\MapArithmetic{m}$-prime. If~\ref{item:prime_m_cliffs_c} holds, by
    Lemma~\ref{lem:decomposition_zero}, there exists $w'' \in \SetCliff_\MapArithmetic{m}$
    and $i \in \Han{\Brr{w''}}$ such that $w' = w'' \WhiteSquare_i 0$. Therefore, $w = w' a
    = w'' a \WhiteSquare_i 0$. Since $a < (\Length(w) - 2)m + 1$, $w'' a$ is an
    $\MapArithmetic{m}$-cliff, showing that $w$ is not $\MapArithmetic{m}$-prime. We have
    shown that in all possible situations, $w$ is not $\MapArithmetic{m}$-prime,
    establishing the equivalence of the statement of the lemma.
\end{proof}

\begin{Theorem} \label{thm:cardinal_minimal_generating_set_cliff_m}
    For any $m \geq 0$, $\# \GeneratingSet_\MapArithmetic{m}(1) = 0$, $\#
    \GeneratingSet_\MapArithmetic{m}(2) = 1$, and, for any $n \geq 3$,
    \begin{equation} \label{equ:cardinal_minimal_generating_set_cliff_m}
        \# \GeneratingSet_\MapArithmetic{m}(n)
        = \frac{m}{m + 1} \Par{\# \SetCliff_\MapArithmetic{m}(n)}.
    \end{equation}
\end{Theorem}
\begin{proof}
    First, since $\GeneratingSet_\MapArithmetic{m}(1) = \emptyset$ and
    $\GeneratingSet_\MapArithmetic{m}(2) = \Bra{\BasisE_0}$, the first part of the statement
    holds. By using the recursive description of $\MapArithmetic{m}$-prime
    $\MapArithmetic{m}$-cliffs provided by Lemma~\ref{lem:prime_m_cliffs}, for any $n \geq
    3$,
    \begin{equation}
        \# \SetPrime_\MapArithmetic{m}(n) =
        m \Par{\# \SetCliff_\MapArithmetic{m}(n - 1)
            - \# \SetPrime_\MapArithmetic{m}(n - 1)}
        + m(n - 2) \# \SetPrime_\MapArithmetic{m}(n - 1).
    \end{equation}
    Therefore,
    \begin{equation} \label{equ:cardinal_minimal_generating_set_cliff_m_1}
        \# \SetPrime_\MapArithmetic{m}(n) =
        m \Par{\# \SetCliff_\MapArithmetic{m}(n - 1)
            + (n - 3) \# \SetPrime_\MapArithmetic{m}(n - 1)}.
    \end{equation}
    By induction on $n \geq 3$, we obtain
    from~\eqref{equ:cardinal_minimal_generating_set_cliff_m_1} that
    \begin{math}
        \# \SetPrime_\MapArithmetic{m}(n)
        = \frac{m}{m + 1} \Par{\# \SetCliff_\MapArithmetic{m}(n)}.
    \end{math}
    Finally, since by Proposition~\ref{prop:minimal_generating_set_cliffs}, for any $n \geq
    1$, $\# \SetPrime_\MapArithmetic{m}(n) = \# \GeneratingSet_\MapArithmetic{m}(n)$,
    \eqref{equ:cardinal_minimal_generating_set_cliff_m} follows.
\end{proof}

By Theorem~\ref{thm:cardinal_minimal_generating_set_cliff_m}, the sequences of the
numbers of elements of $\GeneratingSet_\MapArithmetic{m}$, $m \in \HanL{2}$, counted w.r.t.\
their arities start with
\begin{subequations}
\begin{equation}
    0, 1, 0, 0, 0, 0, 0, 0, 0, \qquad m = 0,
\end{equation}
\begin{equation}
    0, 1, 1, 3, 12, 60, 360, 2520, 20160,
    \qquad m = 1,
\end{equation}
\begin{equation}
    0, 1, 2, 10, 70, 630, 6930, 90090, 1351350,
    \qquad m = 2.
\end{equation}
\end{subequations}
Observe that the minimal generating set $\GeneratingSet_\MapArithmetic{1}$ of
$\SpaceCliff_\MapArithmetic{1}$ is in one-to-one correspondence with the set of even
permutations. The second and third sequences are, respectively, Sequences~\OEIS{A001710}
and~\OEIS{A293962} of~\cite{Slo}.

Here are the list of the first elements of generating families of the relation spaces
$\RelationSpace_\MapArithmetic{m}$, $m \in \HanL{2}$, up to arity~$4$:
\begin{subequations}
\begin{equation}
    \BasisE_0 \Composition_1 \BasisE_0 - \BasisE_0 \Composition_2 \BasisE_0,
    \qquad m = 0,
\end{equation}
\begin{equation}
    \BasisE_0 \Composition_1 \BasisE_0 - \BasisE_0 \Composition_2 \BasisE_0,
    \quad
    \BasisE_{01} \Composition_1 \BasisE_0 - \BasisE_0 \Composition_2 \BasisE_{01},
    \enspace
    \BasisE_0 \Composition_2 \BasisE_{01} - \BasisE_{01} \Composition_2 \BasisE_0,
    \enspace
    \BasisE_0 \Composition_1 \BasisE_{01} - \BasisE_{01} \Composition_3 \BasisE_0,
    \qquad m = 1,
\end{equation}
\begin{multline}
    \BasisE_0 \Composition_1 \BasisE_0 - \BasisE_0 \Composition_2 \BasisE_0,
    \\
    \BasisE_{02} \Composition_1 \BasisE_0 - \BasisE_0 \Composition_2 \BasisE_{02},
    \enspace
    \BasisE_{02} \Composition_2 \BasisE_0 - \BasisE_0 \Composition_2 \BasisE_{02},
    \enspace
    \BasisE_0 \Composition_1 \BasisE_{02} - \BasisE_{02} \Composition_3 \BasisE_0,
    \\
    \BasisE_{01} \Composition_1 \BasisE_0 - \BasisE_0 \Composition_2 \BasisE_{01},
    \enspace
    \BasisE_{01} \Composition_2 \BasisE_0 - \BasisE_0 \Composition_2 \BasisE_{01},
    \enspace
    \BasisE_0 \Composition_1 \BasisE_{01} - \BasisE_{01} \Composition_3 \BasisE_0,
    \qquad m = 2.
\end{multline}
\end{subequations}
The space $\RelationSpace_\MapArithmetic{0}$ is finitely generated while, by
Proposition~\ref{prop:not_finiteness_relation_space_cliff_operads},
$\RelationSpace_\MapArithmetic{1}$ and $\RelationSpace_\MapArithmetic{2}$ are not. Despite
what these lists of nontrivial relations suggest, for any $m \geq 1$,
$\SpaceCliff_\MapArithmetic{m}$ is not a quadratic operad. Indeed, for any $m \geq 1$,
$\RelationSpace_\MapArithmetic{m}$ contains the nontrivial relation
\begin{equation}
    \BasisE_{002} \Composition_3 \BasisE_{01}
    -
    \Par{\BasisE_0 \Composition_2 \BasisE_0} \Composition_3 \BasisE_{012}
\end{equation}
of arity $6$ which is nonhomogeneous in terms of degrees and nonquadratic. These spaces
$\RelationSpace_\MapArithmetic{m}$, $m \geq 1$, seem hard to describe. With the help of the
computer, we obtain that the sequences of the dimensions of
$\RelationSpace_\MapArithmetic{m}$, $m \in \HanL{2}$, begin by
\begin{subequations}
\begin{equation}
    0, 0, 1, 0, 0, 0, 0,
    \qquad m = 0,
\end{equation}
\begin{equation}
    0, 0, 1, 3, 13, 65, 372, 2424,
    \qquad m = 1,
\end{equation}
\begin{equation}
    0, 0, 1, 6, 44, 378, 3788,
    \qquad m = 2.
\end{equation}
\end{subequations}
For the time being, the last two sequences do not appear in~\cite{Slo}.

\subsection{Operads on hills}
We provide here some properties about quotient operads of $\SpaceCliff_\delta$ whose bases
are indexed by $\delta$-hills.

\subsubsection{General properties}
Let us begin by presenting some general properties of $\delta$-hills and of the quotient
operads $\SpaceCliff_\SubFamilly$ where $\SubFamilly$ is a set of $\delta$-hills and
$\delta$ is a weakly increasing range map.

\begin{Lemma} \label{lem:closures_hill}
    For any weakly increasing range map $\delta$, $\SetHill_\delta$ is closed by subword
    reduction.
\end{Lemma}
\begin{proof}
    This is a straightforward consequence of the fact that any subword $w'$ of a
    $\delta$-hill $w$ is weakly increasing and of the fact that since $\delta$ is weakly
    increasing, $\Reduction_\delta\Par{w'}$ remains weakly increasing.
\end{proof}

By Proposition~\ref{prop:quotient_cliff_operads} and Lemma~\ref{lem:closures_hill},
$\SpaceHill_\delta := \SpaceCliff_{\SetHill_\delta}$ is an operad. For instance, in
$\SpaceHill_{1334^\omega}$,
\begin{subequations}
\begin{equation}
    \BasisE_{0234} \Composition_2 \BasisE_{112} = \BasisE_{0112234},
\end{equation}
\begin{equation}
    \BasisF_{0234} \Composition_2 \BasisF_{112} = \BasisF_{0112234} + \BasisF_{0112244},
\end{equation}
\begin{equation}
    \BasisH_{0234} \Composition_2 \BasisH_{112} = \BasisH_{0112244},
\end{equation}
\end{subequations}
\begin{subequations}
\begin{equation}
    \BasisE_{0234} \Composition_3 \BasisE_{112} = \BasisE_{0222234},
\end{equation}
\begin{equation}
    \BasisF_{0234} \Composition_3 \BasisF_{112} = 0,
\end{equation}
\begin{equation}
    \BasisH_{0234} \Composition_3 \BasisH_{112} = \BasisH_{0111244}.
\end{equation}
\end{subequations}

By Theorem~\ref{thm:composition_f_basis_quotient}, the support of any partial composition
over the $\BasisF$-basis in this operad is an interval of the $\delta$-hill poset
$\Par{\SetHill_\delta(n), \Leq}$, $n \geq 1$, introduced in~\cite{CG20,CG22}. As shown
here, this poset is also a sublattice of $\Par{\SetCliff_\delta(n), \Leq}$. In order to
express the partial composition of $\SpaceHill_\delta$ over the $\BasisE$-basis and the
$\BasisH$-basis, let us introduce the following notations. For any $w \in \SetCliff_\delta$,
let $\overline{w}$ (resp.\ $\underline{w}$) be the $\delta$-hill defined for any $j \in
[\Length(w)]$ by $\overline{w}(j) := \max \Bra{w(1), \dots, w(j)}$ (resp.\ $\underline{w}(j)
:= \min \Bra{w(j), \dots, w(\Length(w))}$).

\begin{Proposition} \label{prop:compositions_e_h_bases_hill_weakly_increasing}
    Let $\delta$ be a weakly increasing range map. For any $u, v \in \SetHill_\delta$ and $i
    \in [|u|]$,
    \begin{equation}
        \BasisE_u \Composition_i \BasisE_v = \BasisE_{\overline{u \WhiteSquare_i v}}
        \quad \mbox{and} \quad
        \BasisH_u \Composition_i \BasisH_v = \BasisH_{\underline{u \BlackSquare_i v}}.
    \end{equation}
\end{Proposition}
\begin{proof}
    For any $w \in \SetCliff_\delta$, let $X_w := \Bra{w' \in \SetHill_\delta : w \Leq w'}$.
    Observe that since $X_w$ contains $\delta(1, \Length(w))$, $X_w \ne \emptyset$.
    Moreover, since $\delta$ is weakly increasing, $u \WhiteSquare_i v$ is a $\delta$-cliff.
    For these reasons, and due to the fact that $\SetCliff_\delta(n)$, $n \geq 1$, is a
    sublattice of $\SetHill_\delta(n)$, by
    Proposition~\ref{prop:composition_e_basis_quotient},
    \begin{equation}
        \BasisE_u \Composition_i \BasisE_v
        = \BasisE_{\Meet_{\SetHill_\delta}\Par{X_{u \WhiteSquare_i v}}}.
    \end{equation}
    The stated formula for the partial composition of two elements expressed over the
    $\BasisE$-basis of $\SpaceHill_\delta$ is now the consequence of the fact that for any
    $w \in \SetCliff_\delta$, $\Meet_{\SetHill_\delta}\Par{X_w} = \overline{w}$, which
    follows by induction on the length of~$w$.

    The formula for the partial composition of two elements expressed over the
    $\BasisH$-basis follows from similar arguments.
\end{proof}

Proposition~\ref{prop:compositions_e_h_bases_hill_weakly_increasing} implies that
$\SpaceHill_\delta$ is a set-operad.

\subsubsection{On constant range maps}
Let us study the operads $\SpaceHill_\MapConstant{c}$, $c \geq 0$. The map $\ToPath$ (see
Section~\ref{subsubsec:rectangular_paths}) allows us to interpret any $\MapConstant{c}$-hill
as a $c$-rectangular path. Therefore, $\SpaceHill_{\MapConstant{c}}$ can be seen as an
operad on $c$-rectangular paths. For instance, in $\SpaceHill_\MapConstant{2}$,
\begin{subequations}
\begin{equation}
    \BasisE_{
    \scalebox{.7}{
    \begin{tikzpicture}[Centering,scale=.35]
        \draw[PathGrid](0,0)grid(3,2);
        \node[PathNode](0)at(0,0){};
        \node[PathNode](1)at(1,0){};
        \node[PathNode](2)at(1,1){};
        \node[PathNode](3)at(1,2){};
        \node[PathNode](4)at(2,2){};
        \node[PathNode](5)at(3,2){};
        \draw[PathStep](0)--(1);
        \draw[PathStep](1)--(2);
        \draw[PathStep](2)--(3);
        \draw[PathStep](3)--(4);
        \draw[PathStep](4)--(5);
    \end{tikzpicture}}}
    \Composition_2
    \BasisE_{
    \scalebox{.7}{
    \begin{tikzpicture}[Centering,scale=.35]
        \draw[PathGrid](0,0)grid(2,2);
        \node[PathNode](0)at(0,0){};
        \node[PathNode](1)at(0,1){};
        \node[PathNode](2)at(1,1){};
        \node[PathNode](3)at(1,2){};
        \node[PathNode](4)at(2,2){};
        \draw[PathStep](0)--(1);
        \draw[PathStep](1)--(2);
        \draw[PathStep](2)--(3);
        \draw[PathStep](3)--(4);
    \end{tikzpicture}}}
    =
    \BasisE_{
    \scalebox{.7}{
    \begin{tikzpicture}[Centering,scale=.35]
        \draw[PathGrid](0,0)grid(5,2);
        \node[PathNode](0)at(0,0){};
        \node[PathNode](1)at(1,0){};
        \node[PathNode](2)at(1,1){};
        \node[PathNode](3)at(2,1){};
        \node[PathNode](4)at(2,2){};
        \node[PathNode](5)at(3,2){};
        \node[PathNode](6)at(4,2){};
        \node[PathNode](7)at(5,2){};
        \draw[PathStep](0)--(1);
        \draw[PathStep](1)--(2);
        \draw[PathStep](2)--(3);
        \draw[PathStep](3)--(4);
        \draw[PathStep](4)--(5);
        \draw[PathStep](5)--(6);
        \draw[PathStep](6)--(7);
    \end{tikzpicture}}},
\end{equation}
\begin{equation}
    \BasisF_{
    \scalebox{.7}{
    \begin{tikzpicture}[Centering,scale=.35]
        \draw[PathGrid](0,0)grid(3,2);
        \node[PathNode](0)at(0,0){};
        \node[PathNode](1)at(1,0){};
        \node[PathNode](2)at(1,1){};
        \node[PathNode](3)at(1,2){};
        \node[PathNode](4)at(2,2){};
        \node[PathNode](5)at(3,2){};
        \draw[PathStep](0)--(1);
        \draw[PathStep](1)--(2);
        \draw[PathStep](2)--(3);
        \draw[PathStep](3)--(4);
        \draw[PathStep](4)--(5);
    \end{tikzpicture}}}
    \Composition_2
    \BasisF_{
    \scalebox{.7}{
    \begin{tikzpicture}[Centering,scale=.35]
        \draw[PathGrid](0,0)grid(2,2);
        \node[PathNode](0)at(0,0){};
        \node[PathNode](1)at(0,1){};
        \node[PathNode](2)at(1,1){};
        \node[PathNode](3)at(1,2){};
        \node[PathNode](4)at(2,2){};
        \draw[PathStep](0)--(1);
        \draw[PathStep](1)--(2);
        \draw[PathStep](2)--(3);
        \draw[PathStep](3)--(4);
    \end{tikzpicture}}}
    =
    \BasisF_{
    \scalebox{.7}{
    \begin{tikzpicture}[Centering,scale=.35]
        \draw[PathGrid](0,0)grid(5,2);
        \node[PathNode](0)at(0,0){};
        \node[PathNode](1)at(1,0){};
        \node[PathNode](2)at(1,1){};
        \node[PathNode](3)at(2,1){};
        \node[PathNode](4)at(2,2){};
        \node[PathNode](5)at(3,2){};
        \node[PathNode](6)at(4,2){};
        \node[PathNode](7)at(5,2){};
        \draw[PathStep](0)--(1);
        \draw[PathStep](1)--(2);
        \draw[PathStep](2)--(3);
        \draw[PathStep](3)--(4);
        \draw[PathStep](4)--(5);
        \draw[PathStep](5)--(6);
        \draw[PathStep](6)--(7);
    \end{tikzpicture}}},
\end{equation}
\begin{equation}
    \BasisH_{
    \scalebox{.7}{
    \begin{tikzpicture}[Centering,scale=.35]
        \draw[PathGrid](0,0)grid(3,2);
        \node[PathNode](0)at(0,0){};
        \node[PathNode](1)at(1,0){};
        \node[PathNode](2)at(1,1){};
        \node[PathNode](3)at(1,2){};
        \node[PathNode](4)at(2,2){};
        \node[PathNode](5)at(3,2){};
        \draw[PathStep](0)--(1);
        \draw[PathStep](1)--(2);
        \draw[PathStep](2)--(3);
        \draw[PathStep](3)--(4);
        \draw[PathStep](4)--(5);
    \end{tikzpicture}}}
    \Composition_2
    \BasisH_{
    \scalebox{.7}{
    \begin{tikzpicture}[Centering,scale=.35]
        \draw[PathGrid](0,0)grid(2,2);
        \node[PathNode](0)at(0,0){};
        \node[PathNode](1)at(0,1){};
        \node[PathNode](2)at(1,1){};
        \node[PathNode](3)at(1,2){};
        \node[PathNode](4)at(2,2){};
        \draw[PathStep](0)--(1);
        \draw[PathStep](1)--(2);
        \draw[PathStep](2)--(3);
        \draw[PathStep](3)--(4);
    \end{tikzpicture}}}
    =
    \BasisH_{
    \scalebox{.7}{
    \begin{tikzpicture}[Centering,scale=.35]
        \draw[PathGrid](0,0)grid(5,2);
        \node[PathNode](0)at(0,0){};
        \node[PathNode](1)at(1,0){};
        \node[PathNode](2)at(1,1){};
        \node[PathNode](3)at(2,1){};
        \node[PathNode](4)at(2,2){};
        \node[PathNode](5)at(3,2){};
        \node[PathNode](6)at(4,2){};
        \node[PathNode](7)at(5,2){};
        \draw[PathStep](0)--(1);
        \draw[PathStep](1)--(2);
        \draw[PathStep](2)--(3);
        \draw[PathStep](3)--(4);
        \draw[PathStep](4)--(5);
        \draw[PathStep](5)--(6);
        \draw[PathStep](6)--(7);
    \end{tikzpicture}}},
\end{equation}
\end{subequations}
\begin{subequations}
\begin{equation}
    \BasisE_{
    \scalebox{.7}{
    \begin{tikzpicture}[Centering,scale=.35]
        \draw[PathGrid](0,0)grid(5,2);
        \node[PathNode](0)at(0,0){};
        \node[PathNode](1)at(1,0){};
        \node[PathNode](2)at(1,1){};
        \node[PathNode](3)at(2,1){};
        \node[PathNode](4)at(3,1){};
        \node[PathNode](5)at(3,2){};
        \node[PathNode](6)at(3,2){};
        \node[PathNode](7)at(4,2){};
        \node[PathNode](8)at(5,2){};
        \draw[PathStep](0)--(1);
        \draw[PathStep](1)--(2);
        \draw[PathStep](2)--(3);
        \draw[PathStep](3)--(4);
        \draw[PathStep](4)--(5);
        \draw[PathStep](5)--(6);
        \draw[PathStep](6)--(7);
        \draw[PathStep](7)--(8);
    \end{tikzpicture}}}
    \Composition_3
    \BasisE_{
    \scalebox{.7}{
    \begin{tikzpicture}[Centering,scale=.35]
        \draw[PathGrid](0,0)grid(2,2);
        \node[PathNode](0)at(0,0){};
        \node[PathNode](1)at(1,0){};
        \node[PathNode](2)at(2,0){};
        \draw[PathStep](0)--(1);
        \draw[PathStep](1)--(2);
    \end{tikzpicture}}}
    =
    \BasisE_{
    \scalebox{.7}{
    \begin{tikzpicture}[Centering,scale=.35]
        \draw[PathGrid](0,0)grid(7,2);
        \node[PathNode](0)at(0,0){};
        \node[PathNode](1)at(1,0){};
        \node[PathNode](2)at(1,1){};
        \node[PathNode](3)at(2,1){};
        \node[PathNode](4)at(3,1){};
        \node[PathNode](5)at(4,1){};
        \node[PathNode](6)at(5,1){};
        \node[PathNode](7)at(5,2){};
        \node[PathNode](8)at(6,2){};
        \node[PathNode](9)at(7,2){};
        \draw[PathStep](0)--(1);
        \draw[PathStep](1)--(2);
        \draw[PathStep](2)--(3);
        \draw[PathStep](3)--(4);
        \draw[PathStep](4)--(5);
        \draw[PathStep](5)--(6);
        \draw[PathStep](6)--(7);
        \draw[PathStep](7)--(8);
        \draw[PathStep](8)--(9);
    \end{tikzpicture}}},
\end{equation}
\begin{equation}
    \BasisF_{
    \scalebox{.7}{
    \begin{tikzpicture}[Centering,scale=.35]
        \draw[PathGrid](0,0)grid(5,2);
        \node[PathNode](0)at(0,0){};
        \node[PathNode](1)at(1,0){};
        \node[PathNode](2)at(1,1){};
        \node[PathNode](3)at(2,1){};
        \node[PathNode](4)at(3,1){};
        \node[PathNode](5)at(3,2){};
        \node[PathNode](6)at(3,2){};
        \node[PathNode](7)at(4,2){};
        \node[PathNode](8)at(5,2){};
        \draw[PathStep](0)--(1);
        \draw[PathStep](1)--(2);
        \draw[PathStep](2)--(3);
        \draw[PathStep](3)--(4);
        \draw[PathStep](4)--(5);
        \draw[PathStep](5)--(6);
        \draw[PathStep](6)--(7);
        \draw[PathStep](7)--(8);
    \end{tikzpicture}}}
    \Composition_3
    \BasisF_{
    \scalebox{.7}{
    \begin{tikzpicture}[Centering,scale=.35]
        \draw[PathGrid](0,0)grid(2,2);
        \node[PathNode](0)at(0,0){};
        \node[PathNode](1)at(1,0){};
        \node[PathNode](2)at(2,0){};
        \draw[PathStep](0)--(1);
        \draw[PathStep](1)--(2);
    \end{tikzpicture}}}
    =
    0,
\end{equation}
\begin{equation}
    \BasisH_{
    \scalebox{.7}{
    \begin{tikzpicture}[Centering,scale=.35]
        \draw[PathGrid](0,0)grid(5,2);
        \node[PathNode](0)at(0,0){};
        \node[PathNode](1)at(1,0){};
        \node[PathNode](2)at(1,1){};
        \node[PathNode](3)at(2,1){};
        \node[PathNode](4)at(3,1){};
        \node[PathNode](5)at(3,2){};
        \node[PathNode](6)at(3,2){};
        \node[PathNode](7)at(4,2){};
        \node[PathNode](8)at(5,2){};
        \draw[PathStep](0)--(1);
        \draw[PathStep](1)--(2);
        \draw[PathStep](2)--(3);
        \draw[PathStep](3)--(4);
        \draw[PathStep](4)--(5);
        \draw[PathStep](5)--(6);
        \draw[PathStep](6)--(7);
        \draw[PathStep](7)--(8);
    \end{tikzpicture}}}
    \Composition_3
    \BasisH_{
    \scalebox{.7}{
    \begin{tikzpicture}[Centering,scale=.35]
        \draw[PathGrid](0,0)grid(2,2);
        \node[PathNode](0)at(0,0){};
        \node[PathNode](1)at(1,0){};
        \node[PathNode](2)at(2,0){};
        \draw[PathStep](0)--(1);
        \draw[PathStep](1)--(2);
    \end{tikzpicture}}}
    =
    \BasisH_{
    \scalebox{.7}{
    \begin{tikzpicture}[Centering,scale=.35]
        \draw[PathGrid](0,0)grid(7,2);
        \node[PathNode](0)at(0,0){};
        \node[PathNode](1)at(1,0){};
        \node[PathNode](2)at(2,0){};
        \node[PathNode](3)at(3,0){};
        \node[PathNode](4)at(4,0){};
        \node[PathNode](5)at(4,1){};
        \node[PathNode](6)at(5,1){};
        \node[PathNode](7)at(5,2){};
        \node[PathNode](8)at(6,2){};
        \node[PathNode](9)at(7,2){};
        \draw[PathStep](0)--(1);
        \draw[PathStep](1)--(2);
        \draw[PathStep](2)--(3);
        \draw[PathStep](3)--(4);
        \draw[PathStep](4)--(5);
        \draw[PathStep](5)--(6);
        \draw[PathStep](6)--(7);
        \draw[PathStep](7)--(8);
        \draw[PathStep](8)--(9);
    \end{tikzpicture}}}.
\end{equation}
\end{subequations}

The next result provides a finite presentation by generators and relations
of~$\SpaceHill_\MapConstant{c}$.

\begin{Theorem} \label{thm:presentation_hill_constant_operad}
    For any $c \geq 0$, the set $\Bra{\BasisE_0, \dots, \BasisE_c}$ is a minimal generating
    set of the operad $\SpaceHill_\MapConstant{c}$. The space of the nontrivial relations of
    $\SpaceHill_\MapConstant{c}$ is generated by
    \begin{subequations}
    \begin{equation} \label{equ:presentation_hill_constant_operad_1}
        \BasisE_a \Composition_1 \BasisE_b - \BasisE_b \Composition_2 \BasisE_{a'},
        \qquad
        b \in \HanL{c}, \enspace a, a' \in \HanL{b},
    \end{equation}
    \begin{equation} \label{equ:presentation_hill_constant_operad_2}
        \BasisE_b \Composition_1 \BasisE_a - \BasisE_a \Composition_2 \BasisE_b,
        \qquad
        b \in \HanL{c}, \enspace a \in \HanL{b - 1}.
    \end{equation}
    \end{subequations}
    Moreover, $\SpaceHill_\MapConstant{c}$ is a Koszul operad.
\end{Theorem}
\begin{proof}
    By Proposition~\ref{prop:compositions_e_h_bases_hill_weakly_increasing}, for any $w a
    \in \SetHill_\MapConstant{c}$ such that $w \in \SetHill_\MapConstant{c}$ and $a \in
    \HanL{c}$, one has $\BasisE_w \Composition_{|w|} \BasisE_a = \BasisE_{wa}$. Therefore,
    it follows by induction on the arity that the stated set is a generating family of
    $\SpaceHill_\MapConstant{c}$. Its minimality follows from the fact that no $\BasisE_a$,
    $a \in \HanL{c}$, can be written as a partial composition of other elements of this
    family.

    Let us now show that the space of the nontrivial relations of
    $\SpaceHill_\MapConstant{c}$ is generated
    by~\eqref{equ:presentation_hill_constant_operad_1}
    and~\eqref{equ:presentation_hill_constant_operad_2}. Since the evaluations in
    $\SpaceHill_\MapConstant{c}$ of these two families of expressions is $0$, they belong to
    the space of the nontrivial relations of $\SpaceHill_\MapConstant{c}$. Let $\Operad$ be
    the quotient of the free operad generated by the stated minimal generating set by the
    operad ideal generated by~\eqref{equ:presentation_hill_constant_operad_1}
    and~\eqref{equ:presentation_hill_constant_operad_2}. From these relations, it appears
    that a basis of $\Operad$ is formed by the set $T := \Bra{\TreeT_{\Par{\alpha_0, \dots,
    \alpha_c}} : \alpha_i \in \N, i \in \HanL{c}}$ of expressions defined by
    \begin{equation} \label{equ:presentation_hill_constant_operad_3}
        \TreeT_{\Par{\alpha_0, \dots, \alpha_c}} :=
        \underbrace{\BasisE_c \Composition_1 \dots \Composition_1 \BasisE_c}
            _{\alpha_c \mbox{ terms}}
        \Composition_1
        \underbrace{\BasisE_{c - 1} \Composition_1 \dots \Composition_1 \BasisE_{c - 1}}
            _{\alpha_{c - 1} \mbox{ terms}}
        \Composition_1 \dots \Composition_1
        \underbrace{\BasisE_0 \Composition_1 \dots \Composition_1 \BasisE_0}
            _{\alpha_0 \mbox{ terms}}.
    \end{equation}
    For any $n \geq 1$, the bases of $\SpaceHill_\MapConstant{c}$ and of $\Operad$
    restrained to arity $n$ are in one-to-one correspondence: each $\MapConstant{c}$-hill
    $w$ is in correspondence with the sequence $\Par{\alpha_0, \dots, \alpha_c}$ such that
    $w = 0^{\alpha_0} \dots c^{\alpha_c}$. Therefore, the operads $\Operad$ and
    $\SpaceHill_\MapConstant{c}$ are isomorphic and the stated property holds.

    Finally, the set $T$ forms a Poincaré-Birkhoff-Witt basis of
    $\SpaceHill_\MapConstant{c}$, which implies that this operad is Koszul~\cite{Hof10}.
\end{proof}

By Theorem~\ref{thm:presentation_hill_constant_operad}, for any $c \geq 0$ and any $b \in
\HanL{c}$,
\begin{math}
    \BasisE_b \Composition_1 \BasisE_b = \BasisE_b \Composition_2 \BasisE_b,
\end{math}
so that any algebra on the operad $\SpaceHill_\MapConstant{c}$ has $c + 1$ associative
binary products. This property is shared with some algebras appearing in~\cite{Gir16,Gir16b}
(multiassociative algebras), and in~\cite{ZGG20} (matching associative algebras).

Besides, Theorem~\ref{thm:presentation_hill_constant_operad} shows that
$\SpaceHill_\MapConstant{c}$ is a binary and quadratic operad. Therefore,
$\SpaceHill_\MapConstant{c}$ admits a Koszul dual operad ${\SpaceHill_\MapConstant{c}}^!$.
Let us study this operad now.

\begin{Proposition} \label{prop:presentation_hill_constant_operad_dual}
    For any $c \geq 0$, the operad ${\SpaceHill_\MapConstant{c}}^!$ is isomorphic to the
    quotient of the free operad generated by the set $\Bra{\BasisE^\Dual_0, \dots,
    \BasisE^\Dual_c}$ of symbols of arity $2$ by the operad ideal generated by
    \begin{subequations}
    \begin{equation} \label{equ:presentation_hill_constant_operad_dual_1}
        \sum_{a \in \HanL{b}}
        \BasisE^\Dual_a \Composition_1 \BasisE^\Dual_b
        - \BasisE^\Dual_b \Composition_2 \BasisE^\Dual_a,
        \qquad b \in \HanL{c},
    \end{equation}
    \begin{equation} \label{equ:presentation_hill_constant_operad_dual_2}
        \BasisE^\Dual_b \Composition_1 \BasisE^\Dual_a
        - \BasisE^\Dual_a \Composition_2 \BasisE^\Dual_b,
        \qquad b \in \HanL{c}, \enspace a \in \HanL{b - 1}.
    \end{equation}
    \end{subequations}
\end{Proposition}
\begin{proof}
    This is a straightforward computation based upon the presentation by generators and
    relations of $\SpaceHill_\MapConstant{c}$ provided by
    Theorem~\ref{thm:presentation_hill_constant_operad}. The generating family formed
    by~\eqref{equ:presentation_hill_constant_operad_dual_1}
    and~\eqref{equ:presentation_hill_constant_operad_dual_2} for the space of the nontrivial
    relations of ${\SpaceHill_\MapConstant{c}}^!$ is obtained as the annihilator of the
    linear span of~\eqref{equ:presentation_hill_constant_operad_1}
    and~\eqref{equ:presentation_hill_constant_operad_2} w.r.t.\ to an appropriate linear map
    (see~\cite{Gir18} for instance).
\end{proof}

\begin{Proposition} \label{prop:dimensions_hill_constant_operad_dual}
    For any $c \geq 0$ and any $n \geq 1$,
    \begin{math}
        \dim {\SpaceHill_\MapConstant{c}}^!(n) = \FussCatalan_c(n).
    \end{math}
\end{Proposition}
\begin{proof}
    By Theorem~\ref{thm:presentation_hill_constant_operad}, $\SpaceHill_\MapConstant{c}$ is
    a Koszul operad. Hence, by~\cite{GK94}, its Hilbert series $G(t)$ and the Hilbert series
    $F(t)$ of ${\SpaceHill_\MapConstant{c}}^!$ satisfy $F(-G(-t)) = t = G(-F(-t))$.
    By~\eqref{equ:number_rectangles}, one has
    \begin{equation}
        G(t) = \frac{t}{(1 - t)^{c + 1}}
    \end{equation}
    so that
    \begin{equation}
        G(-F(-t)) = \frac{-F(-t)}{(1 - (-F(-t)))^{c + 1}} = t.
    \end{equation}
    Therefore, $F(t)$ satisfies
    \begin{equation}
        F(t) = t (1 + F(t))^{c + 1}.
    \end{equation}
    This expression for $F(t)$ shows that $F(t)$ is the generating series of the graded set
    of all planar rooted trees such that each internal node has $c + 1$ children, where the
    size is given by the number of internal nodes. Hence, $F(t) = \sum_{n \geq 1}
    \FussCatalan_c(n) t^n$, establishing the statement of the proposition.
\end{proof}

Let us consider for all $b \in \HanL{c}$ the elements
\begin{equation}
    \BasisK^\Dual_b := \sum_{a \in \HanL{b}} \BasisE^\Dual_a
\end{equation}
of the operad~${\SpaceHill_\MapConstant{c}}^!$. By triangularity, $\Bra{\BasisK^\Dual_0,
\dots, \BasisK^\Dual_c}$ is a minimal generating set of~${\SpaceHill_\MapConstant{c}}^!$.

\begin{Proposition} \label{prop:set_presentation_hill_constant_operad_dual}
    For any $c \geq 0$, the operad ${\SpaceHill_\MapConstant{c}}^!$ is isomorphic to the
    quotient of the free operad generated by the set $\Bra{\BasisK^\Dual_0, \dots,
    \BasisK^\Dual_c}$ of symbols of arity $2$ by the operad ideal generated by
    \begin{equation} \label{equ:set_presentation_hill_constant_operad_dual}
        \BasisK^\Dual_b \Composition_1 \BasisK^\Dual_a
        - \BasisK^\Dual_a \Composition_2 \BasisK^\Dual_b,
        \qquad
        b \in \HanL{c}, \enspace a \in \HanL{b}.
    \end{equation}
\end{Proposition}
\begin{proof}
    For any $b \in \HanL{c}$ and $a \in \HanL{b}$, we have
    \begin{equation}
        \BasisK^\Dual_b \Composition_1 \BasisK^\Dual_a
        - \BasisK^\Dual_a \Composition_2 \BasisK^\Dual_b
        =
        \sum_{\substack{
            b' \in \HanL{b} \\
            a' \in \HanL{a}
        }}
        \BasisE^\Dual_{b'} \Composition_1 \BasisE^\Dual_{a'}
        -
        \BasisE^\Dual_{a'} \Composition_2 \BasisE^\Dual_{b'},
    \end{equation}
    implying that~\eqref{equ:set_presentation_hill_constant_operad_dual} expresses as a
    linear combination of~\eqref{equ:presentation_hill_constant_operad_dual_1}
    and~\eqref{equ:presentation_hill_constant_operad_dual_2}. Moreover, the space generated
    by~\eqref{equ:presentation_hill_constant_operad_dual_1}
    and~\eqref{equ:presentation_hill_constant_operad_dual_2} has dimension $\binom{c +
    2}{2}$. We can observe that the space generated
    by~\eqref{equ:set_presentation_hill_constant_operad_dual} has the same dimension.
    Therefore, these two spaces are equal. Finally, since by
    Proposition~\ref{prop:presentation_hill_constant_operad_dual},
    \eqref{equ:presentation_hill_constant_operad_dual_1}
    and~\eqref{equ:presentation_hill_constant_operad_dual_2} form a generating family of the
    nontrivial relations of ${\SpaceHill_\MapConstant{c}}^!$, this is also the case for
    \eqref{equ:set_presentation_hill_constant_operad_dual}. This leads to the stated
    presentation by generators and relations of~${\SpaceHill_\MapConstant{c}}^!$.
\end{proof}

By Proposition~\ref{prop:set_presentation_hill_constant_operad_dual},
${\SpaceHill_\MapConstant{1}}^!$ is the duplicial operad~\cite{BF03} and
${\SpaceHill_\MapConstant{2}}^!$ is the triplicial operad~\cite{Ler11}. For any $c \geq 0$,
we call \Def{$c$-supplicial operad} each operad ${\SpaceHill_\MapConstant{c}}^!$. These
operads are therefore natural generalizations of the duplicial and triplicial operads.

\subsubsection{On arithmetic range maps} \label{subsubsec:fuss_catalan_operads}
Let us study the operads $\SpaceHill_\MapArithmetic{m}$, $m \geq 0$. The map $\ToDyck$ (see
Section~\ref{subsubsec:dyck_paths}) allows us to interpret any $\MapArithmetic{m}$-hill as
an $m$-Dyck path. Therefore, $\SpaceHill_{\MapArithmetic{m}}$ can be seen as an operad on
$m$-Dyck paths. For instance, in $\SpaceHill_\MapArithmetic{1}$ we have
\begin{subequations}
\begin{equation}
    \BasisE_{
    \scalebox{.6}{
    \begin{tikzpicture}[Centering,scale=.3]
        \draw[PathGrid](0,0)grid(10,2);
        \node[PathNode](0)at(0,0){};
        \node[PathNode](1)at(1,1){};
        \node[PathNode](2)at(2,0){};
        \node[PathNode](3)at(3,1){};
        \node[PathNode](4)at(4,0){};
        \node[PathNode](5)at(5,1){};
        \node[PathNode](6)at(6,2){};
        \node[PathNode](7)at(7,1){};
        \node[PathNode](8)at(8,2){};
        \node[PathNode](9)at(9,1){};
        \node[PathNode](10)at(10,0){};
        \draw[PathStep](0)--(1);
        \draw[PathStep](1)--(2);
        \draw[PathStep](2)--(3);
        \draw[PathStep](3)--(4);
        \draw[PathStep](4)--(5);
        \draw[PathStep](5)--(6);
        \draw[PathStep](6)--(7);
        \draw[PathStep](7)--(8);
        \draw[PathStep](8)--(9);
        \draw[PathStep](9)--(10);
    \end{tikzpicture}}}
    \Composition_3
    \BasisE_{
    \scalebox{.6}{
    \begin{tikzpicture}[Centering,scale=.3]
        \draw[PathGrid](0,0)grid(4,1);
        \node[PathNode](0)at(0,0){};
        \node[PathNode](1)at(1,1){};
        \node[PathNode](2)at(2,0){};
        \node[PathNode](3)at(3,1){};
        \node[PathNode](4)at(4,0){};
        \draw[PathStep](0)--(1);
        \draw[PathStep](1)--(2);
        \draw[PathStep](2)--(3);
        \draw[PathStep](3)--(4);
    \end{tikzpicture}}}
    =
    \BasisE_{
    \scalebox{.6}{
    \begin{tikzpicture}[Centering,scale=.3]
        \draw[PathGrid](0,0)grid(14,4);
        \node[PathNode](0)at(0,0){};
        \node[PathNode](1)at(1,1){};
        \node[PathNode](2)at(2,0){};
        \node[PathNode](3)at(3,1){};
        \node[PathNode](4)at(4,2){};
        \node[PathNode](5)at(5,3){};
        \node[PathNode](6)at(6,2){};
        \node[PathNode](7)at(7,3){};
        \node[PathNode](8)at(8,4){};
        \node[PathNode](9)at(9,3){};
        \node[PathNode](10)at(10,4){};
        \node[PathNode](11)at(11,3){};
        \node[PathNode](12)at(12,2){};
        \node[PathNode](13)at(13,1){};
        \node[PathNode](14)at(14,0){};
        \draw[PathStep](0)--(1);
        \draw[PathStep](1)--(2);
        \draw[PathStep](2)--(3);
        \draw[PathStep](3)--(4);
        \draw[PathStep](4)--(5);
        \draw[PathStep](5)--(6);
        \draw[PathStep](6)--(7);
        \draw[PathStep](7)--(8);
        \draw[PathStep](8)--(9);
        \draw[PathStep](9)--(10);
        \draw[PathStep](10)--(11);
        \draw[PathStep](11)--(12);
        \draw[PathStep](12)--(13);
        \draw[PathStep](13)--(14);
    \end{tikzpicture}}},
\end{equation}
\begin{equation}
    \BasisF_{
    \scalebox{.6}{
    \begin{tikzpicture}[Centering,scale=.3]
        \draw[PathGrid](0,0)grid(10,2);
        \node[PathNode](0)at(0,0){};
        \node[PathNode](1)at(1,1){};
        \node[PathNode](2)at(2,0){};
        \node[PathNode](3)at(3,1){};
        \node[PathNode](4)at(4,0){};
        \node[PathNode](5)at(5,1){};
        \node[PathNode](6)at(6,2){};
        \node[PathNode](7)at(7,1){};
        \node[PathNode](8)at(8,2){};
        \node[PathNode](9)at(9,1){};
        \node[PathNode](10)at(10,0){};
        \draw[PathStep](0)--(1);
        \draw[PathStep](1)--(2);
        \draw[PathStep](2)--(3);
        \draw[PathStep](3)--(4);
        \draw[PathStep](4)--(5);
        \draw[PathStep](5)--(6);
        \draw[PathStep](6)--(7);
        \draw[PathStep](7)--(8);
        \draw[PathStep](8)--(9);
        \draw[PathStep](9)--(10);
    \end{tikzpicture}}}
    \Composition_3
    \BasisF_{
    \scalebox{.6}{
    \begin{tikzpicture}[Centering,scale=.3]
        \draw[PathGrid](0,0)grid(4,1);
        \node[PathNode](0)at(0,0){};
        \node[PathNode](1)at(1,1){};
        \node[PathNode](2)at(2,0){};
        \node[PathNode](3)at(3,1){};
        \node[PathNode](4)at(4,0){};
        \draw[PathStep](0)--(1);
        \draw[PathStep](1)--(2);
        \draw[PathStep](2)--(3);
        \draw[PathStep](3)--(4);
    \end{tikzpicture}}}
    =
    \BasisF_{
    \scalebox{.6}{
    \begin{tikzpicture}[Centering,scale=.3]
        \draw[PathGrid](0,0)grid(14,4);
        \node[PathNode](0)at(0,0){};
        \node[PathNode](1)at(1,1){};
        \node[PathNode](2)at(2,0){};
        \node[PathNode](3)at(3,1){};
        \node[PathNode](4)at(4,2){};
        \node[PathNode](5)at(5,3){};
        \node[PathNode](6)at(6,2){};
        \node[PathNode](7)at(7,3){};
        \node[PathNode](8)at(8,4){};
        \node[PathNode](9)at(9,3){};
        \node[PathNode](10)at(10,4){};
        \node[PathNode](11)at(11,3){};
        \node[PathNode](12)at(12,2){};
        \node[PathNode](13)at(13,1){};
        \node[PathNode](14)at(14,0){};
        \draw[PathStep](0)--(1);
        \draw[PathStep](1)--(2);
        \draw[PathStep](2)--(3);
        \draw[PathStep](3)--(4);
        \draw[PathStep](4)--(5);
        \draw[PathStep](5)--(6);
        \draw[PathStep](6)--(7);
        \draw[PathStep](7)--(8);
        \draw[PathStep](8)--(9);
        \draw[PathStep](9)--(10);
        \draw[PathStep](10)--(11);
        \draw[PathStep](11)--(12);
        \draw[PathStep](12)--(13);
        \draw[PathStep](13)--(14);
    \end{tikzpicture}}}
    +
    \BasisF_{
    \scalebox{.6}{
    \begin{tikzpicture}[Centering,scale=.3]
        \draw[PathGrid](0,0)grid(14,4);
        \node[PathNode](0)at(0,0){};
        \node[PathNode](1)at(1,1){};
        \node[PathNode](2)at(2,0){};
        \node[PathNode](3)at(3,1){};
        \node[PathNode](4)at(4,2){};
        \node[PathNode](5)at(5,1){};
        \node[PathNode](6)at(6,2){};
        \node[PathNode](7)at(7,3){};
        \node[PathNode](8)at(8,4){};
        \node[PathNode](9)at(9,3){};
        \node[PathNode](10)at(10,4){};
        \node[PathNode](11)at(11,3){};
        \node[PathNode](12)at(12,2){};
        \node[PathNode](13)at(13,1){};
        \node[PathNode](14)at(14,0){};
        \draw[PathStep](0)--(1);
        \draw[PathStep](1)--(2);
        \draw[PathStep](2)--(3);
        \draw[PathStep](3)--(4);
        \draw[PathStep](4)--(5);
        \draw[PathStep](5)--(6);
        \draw[PathStep](6)--(7);
        \draw[PathStep](7)--(8);
        \draw[PathStep](8)--(9);
        \draw[PathStep](9)--(10);
        \draw[PathStep](10)--(11);
        \draw[PathStep](11)--(12);
        \draw[PathStep](12)--(13);
        \draw[PathStep](13)--(14);
    \end{tikzpicture}}}
    +
    \BasisF_{
    \scalebox{.6}{
    \begin{tikzpicture}[Centering,scale=.3]
        \draw[PathGrid](0,0)grid(14,4);
        \node[PathNode](0)at(0,0){};
        \node[PathNode](1)at(1,1){};
        \node[PathNode](2)at(2,0){};
        \node[PathNode](3)at(3,1){};
        \node[PathNode](4)at(4,0){};
        \node[PathNode](5)at(5,1){};
        \node[PathNode](6)at(6,2){};
        \node[PathNode](7)at(7,3){};
        \node[PathNode](8)at(8,4){};
        \node[PathNode](9)at(9,3){};
        \node[PathNode](10)at(10,4){};
        \node[PathNode](11)at(11,3){};
        \node[PathNode](12)at(12,2){};
        \node[PathNode](13)at(13,1){};
        \node[PathNode](14)at(14,0){};
        \draw[PathStep](0)--(1);
        \draw[PathStep](1)--(2);
        \draw[PathStep](2)--(3);
        \draw[PathStep](3)--(4);
        \draw[PathStep](4)--(5);
        \draw[PathStep](5)--(6);
        \draw[PathStep](6)--(7);
        \draw[PathStep](7)--(8);
        \draw[PathStep](8)--(9);
        \draw[PathStep](9)--(10);
        \draw[PathStep](10)--(11);
        \draw[PathStep](11)--(12);
        \draw[PathStep](12)--(13);
        \draw[PathStep](13)--(14);
    \end{tikzpicture}}},
\end{equation}
\begin{equation}
    \BasisH_{
    \scalebox{.6}{
    \begin{tikzpicture}[Centering,scale=.3]
        \draw[PathGrid](0,0)grid(10,2);
        \node[PathNode](0)at(0,0){};
        \node[PathNode](1)at(1,1){};
        \node[PathNode](2)at(2,0){};
        \node[PathNode](3)at(3,1){};
        \node[PathNode](4)at(4,0){};
        \node[PathNode](5)at(5,1){};
        \node[PathNode](6)at(6,2){};
        \node[PathNode](7)at(7,1){};
        \node[PathNode](8)at(8,2){};
        \node[PathNode](9)at(9,1){};
        \node[PathNode](10)at(10,0){};
        \draw[PathStep](0)--(1);
        \draw[PathStep](1)--(2);
        \draw[PathStep](2)--(3);
        \draw[PathStep](3)--(4);
        \draw[PathStep](4)--(5);
        \draw[PathStep](5)--(6);
        \draw[PathStep](6)--(7);
        \draw[PathStep](7)--(8);
        \draw[PathStep](8)--(9);
        \draw[PathStep](9)--(10);
    \end{tikzpicture}}}
    \Composition_3
    \BasisH_{
    \scalebox{.6}{
    \begin{tikzpicture}[Centering,scale=.3]
        \draw[PathGrid](0,0)grid(4,1);
        \node[PathNode](0)at(0,0){};
        \node[PathNode](1)at(1,1){};
        \node[PathNode](2)at(2,0){};
        \node[PathNode](3)at(3,1){};
        \node[PathNode](4)at(4,0){};
        \draw[PathStep](0)--(1);
        \draw[PathStep](1)--(2);
        \draw[PathStep](2)--(3);
        \draw[PathStep](3)--(4);
    \end{tikzpicture}}}
    =
    \BasisH_{
    \scalebox{.6}{
    \begin{tikzpicture}[Centering,scale=.3]
        \draw[PathGrid](0,0)grid(14,4);
        \node[PathNode](0)at(0,0){};
        \node[PathNode](1)at(1,1){};
        \node[PathNode](2)at(2,0){};
        \node[PathNode](3)at(3,1){};
        \node[PathNode](4)at(4,0){};
        \node[PathNode](5)at(5,1){};
        \node[PathNode](6)at(6,2){};
        \node[PathNode](7)at(7,3){};
        \node[PathNode](8)at(8,4){};
        \node[PathNode](9)at(9,3){};
        \node[PathNode](10)at(10,4){};
        \node[PathNode](11)at(11,3){};
        \node[PathNode](12)at(12,2){};
        \node[PathNode](13)at(13,1){};
        \node[PathNode](14)at(14,0){};
        \draw[PathStep](0)--(1);
        \draw[PathStep](1)--(2);
        \draw[PathStep](2)--(3);
        \draw[PathStep](3)--(4);
        \draw[PathStep](4)--(5);
        \draw[PathStep](5)--(6);
        \draw[PathStep](6)--(7);
        \draw[PathStep](7)--(8);
        \draw[PathStep](8)--(9);
        \draw[PathStep](9)--(10);
        \draw[PathStep](10)--(11);
        \draw[PathStep](11)--(12);
        \draw[PathStep](12)--(13);
        \draw[PathStep](13)--(14);
    \end{tikzpicture}}},
\end{equation}
\end{subequations}
\begin{subequations}
\begin{equation}
    \BasisE_{
    \scalebox{.6}{
    \begin{tikzpicture}[Centering,scale=.3]
        \draw[PathGrid](0,0)grid(10,2);
        \node[PathNode](0)at(0,0){};
        \node[PathNode](1)at(1,1){};
        \node[PathNode](2)at(2,0){};
        \node[PathNode](3)at(3,1){};
        \node[PathNode](4)at(4,2){};
        \node[PathNode](5)at(5,1){};
        \node[PathNode](6)at(6,0){};
        \node[PathNode](7)at(7,1){};
        \node[PathNode](8)at(8,0){};
        \node[PathNode](9)at(9,1){};
        \node[PathNode](10)at(10,0){};
        \draw[PathStep](0)--(1);
        \draw[PathStep](1)--(2);
        \draw[PathStep](2)--(3);
        \draw[PathStep](3)--(4);
        \draw[PathStep](4)--(5);
        \draw[PathStep](5)--(6);
        \draw[PathStep](6)--(7);
        \draw[PathStep](7)--(8);
        \draw[PathStep](8)--(9);
        \draw[PathStep](9)--(10);
    \end{tikzpicture}}}
    \Composition_3
    \BasisE_{
    \scalebox{.6}{
    \begin{tikzpicture}[Centering,scale=.3]
        \draw[PathGrid](0,0)grid(4,2);
        \node[PathNode](0)at(0,0){};
        \node[PathNode](1)at(1,1){};
        \node[PathNode](2)at(2,2){};
        \node[PathNode](3)at(3,1){};
        \node[PathNode](4)at(4,0){};
        \draw[PathStep](0)--(1);
        \draw[PathStep](1)--(2);
        \draw[PathStep](2)--(3);
        \draw[PathStep](3)--(4);
    \end{tikzpicture}}}
    =
    \BasisE_{
    \scalebox{.6}{
    \begin{tikzpicture}[Centering,scale=.3]
        \draw[PathGrid](0,0)grid(14,4);
        \node[PathNode](0)at(0,0){};
        \node[PathNode](1)at(1,1){};
        \node[PathNode](2)at(2,0){};
        \node[PathNode](3)at(3,1){};
        \node[PathNode](4)at(4,2){};
        \node[PathNode](5)at(5,3){};
        \node[PathNode](6)at(6,4){};
        \node[PathNode](7)at(7,3){};
        \node[PathNode](8)at(8,2){};
        \node[PathNode](9)at(9,3){};
        \node[PathNode](10)at(10,2){};
        \node[PathNode](11)at(11,3){};
        \node[PathNode](12)at(12,2){};
        \node[PathNode](13)at(13,1){};
        \node[PathNode](14)at(14,0){};
        \draw[PathStep](0)--(1);
        \draw[PathStep](1)--(2);
        \draw[PathStep](2)--(3);
        \draw[PathStep](3)--(4);
        \draw[PathStep](4)--(5);
        \draw[PathStep](5)--(6);
        \draw[PathStep](6)--(7);
        \draw[PathStep](7)--(8);
        \draw[PathStep](8)--(9);
        \draw[PathStep](9)--(10);
        \draw[PathStep](10)--(11);
        \draw[PathStep](11)--(12);
        \draw[PathStep](12)--(13);
        \draw[PathStep](13)--(14);
    \end{tikzpicture}}},
\end{equation}
\begin{equation}
    \BasisF_{
    \scalebox{.6}{
    \begin{tikzpicture}[Centering,scale=.3]
        \draw[PathGrid](0,0)grid(10,2);
        \node[PathNode](0)at(0,0){};
        \node[PathNode](1)at(1,1){};
        \node[PathNode](2)at(2,0){};
        \node[PathNode](3)at(3,1){};
        \node[PathNode](4)at(4,2){};
        \node[PathNode](5)at(5,1){};
        \node[PathNode](6)at(6,0){};
        \node[PathNode](7)at(7,1){};
        \node[PathNode](8)at(8,0){};
        \node[PathNode](9)at(9,1){};
        \node[PathNode](10)at(10,0){};
        \draw[PathStep](0)--(1);
        \draw[PathStep](1)--(2);
        \draw[PathStep](2)--(3);
        \draw[PathStep](3)--(4);
        \draw[PathStep](4)--(5);
        \draw[PathStep](5)--(6);
        \draw[PathStep](6)--(7);
        \draw[PathStep](7)--(8);
        \draw[PathStep](8)--(9);
        \draw[PathStep](9)--(10);
    \end{tikzpicture}}}
    \Composition_3
    \BasisF_{
    \scalebox{.6}{
    \begin{tikzpicture}[Centering,scale=.3]
        \draw[PathGrid](0,0)grid(4,2);
        \node[PathNode](0)at(0,0){};
        \node[PathNode](1)at(1,1){};
        \node[PathNode](2)at(2,2){};
        \node[PathNode](3)at(3,1){};
        \node[PathNode](4)at(4,0){};
        \draw[PathStep](0)--(1);
        \draw[PathStep](1)--(2);
        \draw[PathStep](2)--(3);
        \draw[PathStep](3)--(4);
    \end{tikzpicture}}}
    =
    0,
\end{equation}
\begin{equation}
    \BasisH_{
    \scalebox{.6}{
    \begin{tikzpicture}[Centering,scale=.3]
        \draw[PathGrid](0,0)grid(10,2);
        \node[PathNode](0)at(0,0){};
        \node[PathNode](1)at(1,1){};
        \node[PathNode](2)at(2,0){};
        \node[PathNode](3)at(3,1){};
        \node[PathNode](4)at(4,2){};
        \node[PathNode](5)at(5,1){};
        \node[PathNode](6)at(6,0){};
        \node[PathNode](7)at(7,1){};
        \node[PathNode](8)at(8,0){};
        \node[PathNode](9)at(9,1){};
        \node[PathNode](10)at(10,0){};
        \draw[PathStep](0)--(1);
        \draw[PathStep](1)--(2);
        \draw[PathStep](2)--(3);
        \draw[PathStep](3)--(4);
        \draw[PathStep](4)--(5);
        \draw[PathStep](5)--(6);
        \draw[PathStep](6)--(7);
        \draw[PathStep](7)--(8);
        \draw[PathStep](8)--(9);
        \draw[PathStep](9)--(10);
    \end{tikzpicture}}}
    \Composition_3
    \BasisH_{
    \scalebox{.6}{
    \begin{tikzpicture}[Centering,scale=.3]
        \draw[PathGrid](0,0)grid(4,2);
        \node[PathNode](0)at(0,0){};
        \node[PathNode](1)at(1,1){};
        \node[PathNode](2)at(2,2){};
        \node[PathNode](3)at(3,1){};
        \node[PathNode](4)at(4,0){};
        \draw[PathStep](0)--(1);
        \draw[PathStep](1)--(2);
        \draw[PathStep](2)--(3);
        \draw[PathStep](3)--(4);
    \end{tikzpicture}}}
    =
    \BasisH_{
    \scalebox{.6}{
    \begin{tikzpicture}[Centering,scale=.3]
        \draw[PathGrid](0,0)grid(14,4);
        \node[PathNode](0)at(0,0){};
        \node[PathNode](1)at(1,1){};
        \node[PathNode](2)at(2,2){};
        \node[PathNode](3)at(3,3){};
        \node[PathNode](4)at(4,4){};
        \node[PathNode](5)at(5,3){};
        \node[PathNode](6)at(6,4){};
        \node[PathNode](7)at(7,3){};
        \node[PathNode](8)at(8,2){};
        \node[PathNode](9)at(9,1){};
        \node[PathNode](10)at(10,0){};
        \node[PathNode](11)at(11,1){};
        \node[PathNode](12)at(12,0){};
        \node[PathNode](13)at(13,1){};
        \node[PathNode](14)at(14,0){};
        \draw[PathStep](0)--(1);
        \draw[PathStep](1)--(2);
        \draw[PathStep](2)--(3);
        \draw[PathStep](3)--(4);
        \draw[PathStep](4)--(5);
        \draw[PathStep](5)--(6);
        \draw[PathStep](6)--(7);
        \draw[PathStep](7)--(8);
        \draw[PathStep](8)--(9);
        \draw[PathStep](9)--(10);
        \draw[PathStep](10)--(11);
        \draw[PathStep](11)--(12);
        \draw[PathStep](12)--(13);
        \draw[PathStep](13)--(14);
    \end{tikzpicture}}}.
\end{equation}
\end{subequations}

Here are the lists of the elements of the minimal generating sets
$\GeneratingSet_{\SetHill_\MapArithmetic{m}}$ of $\SpaceHill_\MapArithmetic{m}$ for $m \in
\HanL{2}$, up to arity $5$:
\begin{subequations}
\begin{equation}
    \BasisE_0, \qquad m = 0,
\end{equation}
\begin{equation}
    \BasisE_0, \quad
    \BasisE_{01}, \quad
    \BasisE_{002}, \BasisE_{012}, \quad
    \BasisE_{0003}, \BasisE_{0013}, \BasisE_{0023}, \BasisE_{0113}, \BasisE_{0123},
    \qquad m = 1,
\end{equation}
\begin{multline}
    \BasisE_0, \quad
    \BasisE_{01}, \BasisE_{02}, \quad
    \BasisE_{003}, \BasisE_{004}, \BasisE_{012}, \BasisE_{013}, \BasisE_{014},
    \BasisE_{023}, \BasisE_{024},
    \\
    \BasisE_{0005}, \BasisE_{0006}, \BasisE_{0015}, \BasisE_{0016}, \BasisE_{0025},
    \BasisE_{0026}, \BasisE_{0034}, \BasisE_{0035}, \BasisE_{0036}, \BasisE_{0045},
    \BasisE_{0046}, \BasisE_{0115}, \BasisE_{0116}, \BasisE_{0123}, \BasisE_{0124},
    \\
    \BasisE_{0125}, \BasisE_{0126}, \BasisE_{0134}, \BasisE_{0135}, \BasisE_{0136},
    \BasisE_{0145}, \BasisE_{0146}, \BasisE_{0225}, \BasisE_{0226}, \BasisE_{0234},
    \BasisE_{0235}, \BasisE_{0236}, \BasisE_{0245}, \BasisE_{0246},
    \qquad m = 2.
\end{multline}
\end{subequations}

\begin{Proposition} \label{prop:minimal_generating_set_hill_m_1}
    One has
    \begin{math}
        \GeneratingSet_{\SetHill_\MapArithmetic{1}}
        = \Bra{\BasisE_w : w \in \SetHill_\MapArithmetic{1}
        \mbox{ and } w(\Length(w)) = \Length(w) - 1}.
    \end{math}
\end{Proposition}
\begin{proof}
    Let us denote by $G$ the set described in the statement of the proposition. Let
    $\BasisE_w \in G$ and assume that there exist $u, v \in \SetHill_\MapArithmetic{1}$ and
    $i \in [|u|]$ such that $\BasisE_w = \BasisE_u \Composition_i \BasisE_v$. By
    Proposition~\ref{prop:compositions_e_h_bases_hill_weakly_increasing}, $w = \overline{u
    \WhiteSquare_i v}$. Therefore, there is in $u$ or in $v$ a letter of value $\Length(w) -
    1$, implying that if it occurs in $u$, then $v = \epsilon$, or if it occurs in $v$ that
    $u = \epsilon$. For this reason, $\BasisE_w$ admits no nontrivial decompositions, so
    that $\BasisE_w \in \GeneratingSet_{\SetHill_\MapArithmetic{1}}$.

    Let us now prove that for any $w \in \SetHill_\MapArithmetic{1}$, $\BasisE_w$ belongs to
    the suboperad of $\SpaceHill_\MapArithmetic{1}$ generated by $G$. We proceed by
    recurrence on $|w|$. If $|w| = 1$, then $w = 0$, and since $\BasisE_0 \in G$, the
    property holds. Otherwise, we have $w = w' a$ where $w' \in \SetHill_\MapArithmetic{1}$
    and $a \in \N$. If $a = \Length(w) - 1$, $\BasisE_w \in G$ so that the property holds.
    Otherwise, let $v$ be the subword of $w'$ of length $a$ obtained by scanning $w'$ from
    right to left and by keeping a letter only if it small enough so that $v$ is a
    $\MapArithmetic{1}$-hill. For instance, for $w := 001222356$, we have $w' = 00122235$,
    $a = 6$, and $v = 012235$. By construction, one has $\BasisE_{va} \in G$. Moreover,
    since $w$ is a $\MapArithmetic{1}$-hill, the letters that are in $w$ which have not been
    selected to form $va$ have values that are necessarily in $va$. For this reason, and due
    to the Proposition~\ref{prop:compositions_e_h_bases_hill_weakly_increasing}, $\BasisE_w$
    can be expressed by partial compositions involving $\BasisE_{va}$ and multiple
    occurrences of $\BasisE_0$. Since $\BasisE_0 \in G$, this shows that $\BasisE_w$ belongs
    to the suboperad of $\SpaceHill_\MapArithmetic{1}$ generated by $G$.
\end{proof}

The sequences of the numbers of elements of $\GeneratingSet_{\SetHill_\MapArithmetic{m}}$,
$m \in \HanL{2}$, counted w.r.t.\ their arities start with
\begin{subequations}
\begin{equation}
    0, 1, 0, 0, 0, 0, 0, 0, 0, \qquad m = 0,
\end{equation}
\begin{equation}
    0, 1, 1, 2, 5, 14, 42, 132, 429, \qquad m = 1
\end{equation}
\begin{equation}
    0, 1, 2, 7, 29, 133, 654, 3383, 18179, \qquad m = 2.
\end{equation}
\end{subequations}
We deduce from Proposition~\ref{prop:minimal_generating_set_hill_m_1} that $\#
\GeneratingSet_{\SetHill_\MapArithmetic{1}}(1) = 0$ and for any $n \geq 2$, $\#
\GeneratingSet_{\SetHill_\MapArithmetic{1}}(n) = \FussCatalan_1(n - 2)$. The sequence of the
cardinalities of $\GeneratingSet_{\SetHill_\MapArithmetic{2}}$ does not appear for the time
being in~\cite{Slo}.

Here are the list of the first elements of the generating families of the relation spaces
$\RelationSpace_{\SetHill_\MapArithmetic{m}}$, $m \in \HanL{2}$, up to arity $4$:
\begin{subequations}
\begin{equation}
    \BasisE_0 \Composition_1 \BasisE_0 - \BasisE_0 \Composition_2 \BasisE_0,
    \qquad m = 0,
\end{equation}
\begin{equation}
    \BasisE_0 \Composition_1 \BasisE_0 - \BasisE_0 \Composition_2 \BasisE_0,
    \quad
    \BasisE_{01} \Composition_1 \BasisE_0 - \BasisE_{01} \Composition_2 \BasisE_0,
    \enspace
    \BasisE_0 \Composition_2 \BasisE_{01} - \BasisE_{01} \Composition_2 \BasisE_0,
    \enspace
    \BasisE_0 \Composition_1 \BasisE_{01} - \BasisE_{01} \Composition_3 \BasisE_0,
    \qquad m = 1,
\end{equation}
\begin{multline}
    \BasisE_0 \Composition_1 \BasisE_0 - \BasisE_0 \Composition_2 \BasisE_0,
    \quad
    \BasisE_{01} \Composition_1 \BasisE_0 - \BasisE_{01} \Composition_2 \BasisE_0,
    \enspace
    \BasisE_0 \Composition_2 \BasisE_{01} - \BasisE_{01} \Composition_2 \BasisE_0,
    \enspace
    \BasisE_0 \Composition_1 \BasisE_{01} - \BasisE_{01} \Composition_3 \BasisE_0,
    \\
    \BasisE_{02} \Composition_1 \BasisE_0 - \BasisE_{02} \Composition_2 \BasisE_0,
    \enspace
    \BasisE_0 \Composition_2 \BasisE_{02} - \BasisE_{02} \Composition_2 \BasisE_0,
    \enspace
    \BasisE_0 \Composition_1 \BasisE_{02} - \BasisE_{02} \Composition_3 \BasisE_0,
    \qquad m = 2.
\end{multline}
\end{subequations}
The space $\RelationSpace_{\SetHill_\MapArithmetic{0}}$ is finitely generated. When $m \geq
1$, due to the description of the partial composition map of $\SpaceHill_\MapArithmetic{m}$
provided by Proposition~\ref{prop:compositions_e_h_bases_hill_weakly_increasing}, it appears
that by setting $w := 0 \dots 0 \ \MapArithmetic{m}(\Length(w))$, $\BasisE_w \in
\GeneratingSet_{\SetHill_\MapArithmetic{m}}$. Moreover, since
\begin{math}
    \BasisE_0 \Composition_1 \BasisE_w = \BasisE_w \Composition_{|w|} \BasisE_0,
\end{math}
\begin{equation}
    \BasisE_0 \Composition_1 \BasisE_w - \BasisE_w \Composition_{|w|} \BasisE_0
\end{equation}
is an element of a minimal generating family of
$\RelationSpace_{\SetHill_\MapArithmetic{m}}$. This shows that
$\RelationSpace_{\SetHill_\MapArithmetic{m}}$ is not finitely generated. Besides, despite
what these lists of nontrivial relations suggest, for any $m \geq 1$,
$\RelationSpace_{\SetHill_\MapArithmetic{m}}$ is not a quadratic operad. Indeed, for any $m
\geq 1$, $\RelationSpace_{\SetHill_\MapArithmetic{m}}$ contains the nontrivial relation
\begin{equation}
    \Par{\BasisE_0 \Composition_1 \BasisE_0} \Composition_1 \BasisE_{01}
    -
    \BasisE_{01} \Composition_3 \BasisE_{01}
\end{equation}
of arity $5$ which is nonhomogeneous in terms of degrees and nonquadratic. These spaces
$\RelationSpace_{\SetHill_\MapArithmetic{m}}$, $m \geq 1$, seem hard to describe. With the
help of the computer, we obtain that the sequences of the dimensions of
$\RelationSpace_{\SetHill_\MapArithmetic{m}}$, $m \in \HanL{2}$, begin by
\begin{subequations}
\begin{equation}
    0, 0, 1, 0, 0, 0, 0,
    \qquad m = 0,
\end{equation}
\begin{equation}
    0, 0, 1, 3, 10, 35, 126, 462,
    \qquad m = 1,
\end{equation}
\begin{equation}
    0, 0, 1, 6, 35, 206, 1231,
    \qquad m = 2.
\end{equation}
\end{subequations}
The second seems Sequence~\OEIS{A001700} of~\cite{Slo} and for the time being, the third
sequence does not appear in~\cite{Slo}.

\section*{Conclusion and open questions}
This work endows various set of combinatorial families (integer compositions, permutations,
$m$-increasing trees, $c$-rectangular paths, $m$-Dyck paths, and more generally
$\delta$-cliffs and $\delta$-hills) with operad structures, all being suboperads or
quotients of the interstice operad on nonnegative integers. The main considered operads of
this work fit into the diagram
\begin{equation}
    \begin{tikzpicture}[Centering]
        \node(IN)at(0,-1){$\IntersticeOperad(\N)$};
        \node(Cl_bar_delta)at(2,0){$\SpaceCliff_{\bar{\delta}}$};
        \node(Cl_delta)at(4,0){$\SpaceCliff_\delta$};
        \node(Hi_delta)at(6,0){$\SpaceHill_\delta$};
        \node(Cl_m)at(3,-1){$\SpaceCliff_\MapArithmetic{m}$};
        \node(Cl_c)at(3,-2){$\SpaceCliff_\MapConstant{c}$};
        \node(Hi_m)at(6,-1){$\SpaceHill_\MapArithmetic{m}$};
        \node(Hi_c)at(6,-2){$\SpaceHill_\MapConstant{c}$};
        \draw[Injection](Cl_bar_delta)--(IN);
        \draw[Surjection](Cl_bar_delta)--(Cl_delta);
        \draw[Surjection](Cl_delta)--(Hi_delta);
        \draw[Injection](Cl_m)--(IN);
        \draw[Injection](Cl_c)--(IN);
        \draw[Surjection](Cl_m)--(Hi_m);
        \draw[Surjection](Cl_c)--(Hi_c);
    \end{tikzpicture}
\end{equation}
of injective ($\tikz{\draw[Injection](0,0)--(0.6,0)}$) and surjective
($\tikz{\draw[Surjection](0,0)--(0.6,0)}$) operad morphisms, where $\delta$ is any unimodal
range map. As shown, even if interstice operads have a very simple algebraic structure, the
operads constructed here have more intricate ones since some of them do not admit finite
minimal generating sets, have an infinite family of nontrivial relations, or have some
nontrivial relations which are nonhomogeneous and nonquadratic.

Here is a list of open questions raised by this research:
\begin{enumerate}[fullwidth]
    \item {\bf (Place of $\SpaceCliff_\delta$ in the world of combinatorial operads)} ---
    This consists in exploring other suboperads and quotients of $\SpaceCliff_\delta$ and
    see if $\SpaceCliff_\delta$ contains as substructures some already known operads.
    Moreover, this axis also consists in searching operad morphisms between
    $\SpaceCliff_\delta$ and other operads. More specifically, one can ask for instance
    about such relations between $\SpaceCliff_\MapArithmetic{m}$ and the operad $\As$.

    \item {\bf (Nontrivial relations of $\SpaceCliff_\delta$)} ---
    Proposition~\ref{prop:not_finiteness_relation_space_cliff_operads} provides a sufficient
    condition for the fact that $\RelationSpace_\delta$ is not finitely generated. We can
    ask about a necessary condition for this property.  Moreover, the description and the
    enumeration arity by arity of a minimal generating set for
    $\RelationSpace_\MapArithmetic{m}$, $m \geq 1$, is open.

    \item {\bf (Presentation of $\SpaceHill_\MapArithmetic{m}$)} --- We know that the
    minimal generating set $\GeneratingSet_{\SetHill_\MapArithmetic{1}}$ of
    $\SpaceHill_\MapArithmetic{1}$ is enumerated by a shifted version of Catalan numbers
    (see Section~\ref{subsubsec:fuss_catalan_operads}). The question of the description and
    the enumeration of the minimal generating set
    $\GeneratingSet_{\SetHill_\MapArithmetic{m}}$ of $\SpaceHill_\MapArithmetic{m}$ for $m
    \geq 2$ is open. The analogous question for a minimal generating set
    $\RelationSpace_{\SetHill_\MapArithmetic{m}}$ is also open for $ m \geq 1$. We
    conjecture in particular that $\RelationSpace_{\SetHill_\MapArithmetic{1}}$ is
    enumerated by Sequence~\OEIS{A001700} of~\cite{Slo}.
\end{enumerate}

\bibliographystyle{alpha}
\bibliography{Bibliography}

\end{document}